\documentclass[12pt]{article}
\usepackage{amsmath,amssymb,latexsym}
\usepackage{fullpage}
\usepackage{theorem}

\begin{document}
\date{}
\title{\bf On Spin $L$-Functions for $GSO_{10}$}
\author{ David Ginzburg\\ School of Mathematical Sciences\\
Sackler Faculty of Exact Sciences\\
Tel-Aviv University, Israel 69978\\
and\\ Joseph  Hundley \\
 Penn State University Mathematics Department\\
University Park, State College PA, 16802, USA } \maketitle
\baselineskip=18pt

\abstract{ In this paper we construct a Rankin-Selberg integral
which represents the $Spin_{10}\times St$ $L$-function attached to
the group $GSO_{10}\times PGL_2$. We use this integral
representation to give some equivalent conditions for a generic
cuspidal representation on $GSO_{10}$ to be a functorial lift from
the group $G_2\times PGL_2$. }

\section{ Introduction }

In this paper we construct a new Rankin-Selberg integral which
represents an $L$ function corresponding to the group
$GSO_{10}\times PGL_2$. More precisely, let $\pi$ denote a generic
cuspidal representation defined on the group $GSO_{10}({\bf A})$,
and let $\tau$ denote a cuspidal representation defined on
$PGL_2({\bf A})$. For simplicity, we shall assume that $\pi$ has a
trivial central character. The $L$ group of $GSO_{10}({\bf
A})\times PGL_2({\bf A})$ is the group $GSpin_{10}({\bf C})\times
SL_2({\bf C})$. Consider the 32 dimensional irreducible
representation of this group given by $Spin_{10}\times St$. Here
$Spin_{10}$ is the 16 dimensional irreducible Spin representation
of $GSpin_{10}({\bf C})$ and $St$ denotes the standard
representation of the group $SL_2({\bf C})$. To this irreducible
representation, one can attach the 32 degree partial $L$ function
denoted by $L^S(Spin_{10}\times St,\pi\times\tau,s)$.

To study this $L$ function we construct a global integral which is
given by
\begin{equation}\label{global0}
\int\limits_{Z({\bf A})GSO_{10}(F)\backslash GSO_{10}({\bf
A})}\varphi_\pi(g) \theta_\tau(g)E(g,s)dg.
\end{equation}
Here $\varphi_\pi$ is a vector in the space of the representation
$\pi$, $Z$ is the center of $GSO_{10}$ and $E(g,s)$ is an
Eisenstein series which described in the beginning of section 2.
The interesting representation in integral \eqref{global0} is the
representation $\theta_\tau$. This representation is constructed
as a residue of an Eisenstein series defined on the group
$GSO_{10}({\bf A})$ as described fully in section 3. The main two
properties of this representation are first, its dependence on the
cuspidal representation $\tau$, and second, its smallness. In
fact, as we prove in section 3, this representation is attached to
the unipotent orbit $(3^31)$.

After showing that this integral is Eulerian, and that we obtain
the above $L$ function, we then study the poles of this $L$
function. We show that it can have at most a simple pole at $s=1$.
It is well known, see \cite{K}, that the stabilizer of a generic
point in the space of the representation $Spin_{10}({\bf C})\times
SL_2({\bf C})$, is the group $G_2({\bf C})\times SL_2({\bf C})$.
Hence one expects that the above $L$ function will have a simple
pole at $s=1$ if and only if the cuspidal representation $\pi$ is
the functorial lift from $G_2\times PGL_2$. In fact this is
exactly what we prove. Indeed, the main result of this paper is
given by the following theorem, which is stated and proved in
section 6. (See section 6 for precise notations).\\
{\bf Main Theorem:} {\sl Let $\pi$ be an irreducible generic
cuspidal representation of the group $GSO_{10}(\bf A)$ which has a
trivial central character. Then the following are equivalent:\\
1) The partial $L$ function $L^S(Spin_{10}\times St,\pi\times
\tau, 2s-1/2)$ has a simple pole at $s=3/4$.\\
2) The period integral
\begin{equation}
\int\limits_{SO_{10}(F)\backslash SO_{10}({\bf A})}\varphi_\pi(g)
\theta_\tau(g)\theta(g)dg
\end{equation}
is nonzero for some choice of data.\\
3) There exists a generic cuspidal representation $\sigma$ of the
exceptional group $G_2({\bf A})$ such that $\pi$ is the weak lift
from the representation $\sigma\times\tau$ of the group $G_2({\bf
A})\times PGL_2({\bf A})$. }

We now describe the content of the paper. In section two we
introduce the global integral we consider, and show that it is
Eulerian. We try to do it in an abstract way. Indeed, we assume
the existence of a representation $\theta$ defined on the group
$GSO_{10}(\bf A)$ which is attached to the unipotent orbit
$(3^31)$. Using that data, we prove in theorem 1, that the global
integral is Eulerian. This method has the advantage that it will
work for any representation which will satisfy the properties
listed in the theorem. This means, that for ${\sl any}$ such
representation the integral will be Eulerian. In section 3 we
construct an example of such a representation, which we denote by
$\theta_\tau$, by means of a residue of an Eisenstein series. We
are well aware of the existence of other such representations as
well.

Section 4 is devoted to the unramified computations.  The local
integral and $L$-function are each expressed as a power series
in $q^{-2s+1/2}$, $\chi$, and $\chi^{-1}$, where $\chi$ is the
character from which $\tau$ is induced.
The coefficients are traces of irreducible representations of
$Spin_{10}({\bf C}).$
The power series are
simplified by multiplying them by a certain polynomial.  The
desired equality is equivalent to a formula for 
tensor products of representations of ``rectangular shape'' which
is due to Okada.

Sections 5 and 6 are devoted to the main theorem. The idea of the
proof is as follows. The hard part in the theorem is to prove that
part 2 implies part 3. To do that we first show that the cuspidal
representation $\pi$ is an endoscopic lifting from a cuspidal
representation $\epsilon\times \tau$ defined on the group
$Sp_6({\bf A})\times SL_2({\bf A})$. This part is not trivial. To
do it we need to use a new construction of a lifting which was
announced in \cite{G2} and is a work in progress in \cite{G3}. In
section 5 we give the precise details we need. We don't give all
the proofs, they will appear in \cite{G3}, but we give enough
details for the reader to get the whole picture. Then, assuming
part two of the above theorem, we show that the representation
$\epsilon$ defined on $Sp_6({\bf A})$ is actually a functorial
lift from a generic representation $\sigma$ defined on the
exceptional group $G_2({\bf A})$.

Finally, we wish to remark that the representation $\theta_\tau$
seems to occur in other constructions of Rankin-Selberg integrals.
Recently, we constructed some other global integrals, using this
representation, which we proved to be Eulerian. It is our full
intention to try to find the $L$ functions they represent.

\section{ A Global Integral}

Let $G$ denote the similitude orthogonal group $GSO_{10}$. Let
$\pi$ denote a generic irreducible cuspidal representation on the
group $G({\bf A})$. We shall assume that it has a trivial central
character. Let $P$ denote the standard maximal parabolic subgroup
of $G$ whose levi part is $GL_1\times GL_5$ which contains the
standard Borel subgroup consisting  of upper unipotent matrices.
We shall denote by $U(P)$ the unipotent radical subgroup of $P$.
Let $E(g,s)$ denote the Eisenstein series defined on the group
$G({\bf A})$ which is associated with the induced representation
$Ind_{P({\bf A})}^{G({\bf A})}\delta_P^s$.

Let $\theta$ denote any automorphic representation defined on the
group $G({\bf A})$. Assume that it has a trivial central
character. For a cusp form $\varphi_\pi$ in the space of $\pi$ we
consider the integral
\begin{equation}\label{global1}
\int\limits_{Z({\bf A})G(F)\backslash G({\bf A})}\varphi_\pi(g)
\theta(g)E(g,s)dg
\end{equation}
Here $Z$ is the center of the group $G$. We are mainly interested
in understanding what conditions we need to impose on the
representation $\theta$ so that integral \eqref{global1} will be
Eulerian with the Whittaker function defined on the representation
$\pi$.

In terms of matrices we consider the group $G$ relative to the
form defined by the matrix $J_{10}$. Here and elsewhere, the
matrix $J_n$ is the $n\times n$ matrix with ones on the other
diagonal and zeros elsewhere. For $1\le i\le 5$, let $\alpha_i$
denote the five simple roots of the group $G$. Let
$x_{\alpha_i}(r)$ denote the one dimensional unipotent subgroup
corresponding to the root $\alpha_i$. We label the roots such that
$$x_{\alpha_1}(r)=I+re'_{1,2}\ \ \ \
x_{\alpha_2}(r)=I+re'_{2,3}\ \ \ \ x_{\alpha_3}(r)=I+re'_{3,4}$$
$$ x_{\alpha_4}(r)=I+re'_{4,5}\ \ \ \ x_{\alpha_5}(r)=I+re'_{4,6}$$
Here $I$ is the $10\times 10$ identity matrix and
$e'_{i,j}=e_{i,j}-e_{11-j,11-i}$. For $1\le i\le 5$ let $w[i]$
denote the simple reflection corresponding to the simple root
$\alpha_i$. We shall write $w[i_1i_2\ldots i_r]$ for
$w[i_1]w[i_2]\ldots w[i_r]$.

Let $\psi$ denote a nontrivial character on the group $F\backslash
{\bf A}$. Let $U$ denote the maximal unipotent subgroup of $G$
which consists of upper triangular matrices. We define two
characters on the group $U$. For an automorphic form $\phi$
defined on $G({\bf A})$ we define its Whittaker model as
$$W_\phi(g)=\int\limits_{U(F)\backslash U({\bf A})}
\phi(ug){\psi'_U}^{-1}(u)du$$ Here, for $u=(u_{i,j})\in U$, we
define $\psi'_U(u)=\psi(u_{1,2}+u_{2,3}+u_{3,4}+u_{4,5}+u_{4,6})$.
Similarly, we define
$$\phi^{U,\psi_U}(g)=\int\limits_{U(F)\backslash U({\bf A})}
\phi(ug)\psi_U(u)du$$ where now
$\psi_U(u)=\psi(u_{1,2}+u_{2,3}+u_{4,5}+u_{4,6})$.

In \cite{G-R-S1} it is explained how to associate with a unipotent
class of a classical group a set of Fourier coefficients. We also
use the notation ${\cal O}_G$ as explained there. Roughly
speaking, if $\sigma$ is an automorphic representation of the
group $G$ one defines ${\cal O}_G(\sigma)$ as follows. It is
defined to be the set of all unipotent classes of $G$ such that
for all ${\cal O'}$ with the property that if ${\cal O'}$ is
greater or not related to a member in ${\cal O}_G(\sigma)$, then
$\sigma$ has no nontrivial Fourier coefficient corresponding to
the unipotent class ${\cal O'}$. Also, the representation $\sigma$
has a nonzero Fourier coefficient corresponding to any unipotent
class ${\cal O}$ in ${\cal O}_G(\sigma)$.

We are now ready to prove the following\\
{\bf Theorem 1:} {\sl We keep the above notations. Let $\theta$ be
an automorphic representation of the group $G({\bf A})$ which
satisfies:\\
1) ${\cal O}_G(\theta)=(3^31)$.\\
2) The integral $\theta^{U,\psi_U}(g)$ is not zero for some choice
of data.\\
Then, integral \eqref{global1} is Eulerian and for $Re(s)$ large
it equals}
\begin{equation}\label{th}
\int\limits_{Z({\bf A})U({\bf A})\backslash G({\bf
A})}\int\limits_{{\bf A}^2} W_\pi(g)\theta^{U,\psi_U}(g)
f_s(w[53]x_{\alpha_3}(r_1)x_{\alpha_3+\alpha_5}(r_2)g)dr_1dr_2dg
\end{equation}
{\sl Proof:} For $Re(s)$ large we unfold integral \eqref{global1}
and we obtain
\begin{equation}\label{global2}
\int\limits_{Z({\bf A})P(F)\backslash G({\bf A})}\varphi_\pi(g)
\theta(g)f_s(g)dg
\end{equation}
Next we expand $\theta(g)$ along the unipotent radical of $P$. The
Levi part of $P$ acts on $U(P)$ with three orbits. Indeed, we may
identify $U(P)$ with all matrices in $G$ of the form
$\begin{pmatrix} I&X\\&I \end{pmatrix}.$ Thus we may parameterize
the orbits by the rank of the matrices $X$, which must be even. It
follows from the cuspidality of $\pi$ that the rank zero and rank
two orbits contribute zero to the integral. Denote
$$\theta^{U(P),\psi}(g)=\int\limits_{X(F)\backslash X({\bf A})}
\theta\left ( \begin{pmatrix} I&X\\&I \end{pmatrix}\right )\psi
(x_{2,1}+x_{3,2})dX.$$ The stabilizer of this character inside the
Levi part of $P$ is the group $GSp_4$. Thus \eqref{global2} equals
\begin{equation}\label{global3}
\int\limits_{Z({\bf A})GSp_4(F)U(P)(F)\backslash G({\bf
A})}\varphi_\pi(g) \theta^{U(P),\psi}(g)f_s(g)dg
\end{equation}
Consider the unipotent group $L$ which is generated by the
matrices $I_{10}+re'_{1,2};\ I_{10}+re'_{1,3};\ I_{10}+re'_{1,4};\
I_{10}+re'_{1,5}$. We expand $\theta^{U(P),\psi}(g)$ along the
group $L(F)\backslash L({\bf A})$. The group $GSp_4$ acts on $L$
with two orbits. By the cuspidality of $\pi$, the trivial orbit
contributes zero to the integral. Denote $V=L\cdot U(P)$, and
define a character $\psi_V$ of $V$ as follows. For $v\in V$ write
$v=lu$ with $l=(l_{i,j})\in L$ and $u=\begin{pmatrix} I&X\\&I
\end{pmatrix}\in U(P)$. We define
$\psi_V(v)=\psi(l_{1,2}+x_{2,1}+x_{3,2}).$ Thus we are left with
the open orbit and hence integral \eqref{global3} equals
\begin{equation}\label{global4}
\int\limits_{Z({\bf A})R(F)V({\bf A})\backslash G({\bf
A})}\varphi_\pi^{V,\psi_V^{-1}}(g) \theta^{V,\psi_V}(g)f_s(g)dg
\end{equation}
where
$$\theta^{V,\psi_V}(g)=\int\limits_{V(F)\backslash V({\bf A})}
\theta(vg)\psi_V(v)dv$$ and similarly we define
$\varphi_\pi^{V,\psi_V^{-1}}(g)$. Also, the group $R=GL_2\cdot Y$
where $$Y=\{
(I_{10}+y_1e'_{2,3})(I_{10}-y_1e'_{4,5})(I_{10}+y_2e'_{2,4})
(I_{10}+y_2e'_{3,5})(I_{10}+y_3e'_{2,5})\}$$ The group $GL_2$ is
embedded in $G$ as all matrices of the form
$diag(|h|,|h|,h,1,|h|,h,1,1)$ where $h\in GL_2$ and $|h|=deth$.

Next we consider the Fourier expansion
$$\theta^{V,\psi_V}(g)=\int\limits_{F\backslash
{\bf
A}}\theta^{V,\psi_V}((I_{10}+y_3e'_{2,5})g)dy_3+\sum_{\alpha\in
F^*}\int\limits_{F\backslash {\bf
A}}\theta^{V,\psi_V}((I_{10}+y_3e'_{2,5})g)\psi(\alpha y_3)dy_3$$
We claim that each summand in the summation over $\alpha$ on the
right hand side is zero. Indeed, since $\alpha\ne 0$ it follows
that the integration over $y_3$ together with the integration over
$V$, produces a Fourier coefficient which corresponds to the
unipotent class $(52^21)$. By assumption 1) in the theorem it
follows that these Fourier coefficients are all zero. Thus we are
left only with the constant term.

Denote $w=w[35].$ Let $Y_c=I_{10}+y_3e'_{2,5}$ and denote
$V_1=wVY_cw^{-1}$. To describe the group $V_1$ in term of
matrices, let $U_1$ denote the standard unipotent radical of the
maximal parabolic subgroup of $G$ whose Levi part is $GL_1\times
GSO_8$. We consider the character $\psi_{U}$ as a character of the
group $U_1$ by restriction.  Next we consider the unipotent
subgroup of $U$ defined by all matrices of the form
$$l(r_1,r_2,r_3,r_4,r_5,r_6,r_7)=I_{10}+r_1e'_{2,3}+r_2e'_{2,4}
+r_3e'_{2,5}+r_4e'_{2,6}+r_5e'_{2,7}+r_6e'_{2,8}+r_7e'_{4,5}$$
Finally, we consider the group of all unipotent matrices of the
form $z(m_1,m_2)=I_{10}+m_1e'_{4,3}+m_2e'_{6,3}$. A matrix
multiplication implies that we have the factorization
$v_1=u_1l(r_1,0,r_3,0,r_5,r_6,r_7)z(m_1,m_2)$ where $v_1\in V_1$
and $u_1\in U_1$. We also consider the character $\psi_U$ as a
character of the group $V_1$ by restriction. Thus, for the above
factorization, we have $\psi_U(v_1)=\psi_U(u_1)\psi(r_1+r_7)$.
From all this we obtain that \eqref{global4} equals
\begin{equation}\label{global5}
\int\limits_{Z({\bf A})GL_2(F)Y(F)Y_c({\bf A})V({\bf A})\backslash
G({\bf A})}\varphi_\pi^{V_1,\psi_U^{-1}}(wg)
\theta^{V_1,\psi_U}(wg)f_s(g)dg
\end{equation}
We proceed with the following Fourier expansion
$$\theta^{V_1,\psi_U}(wg)=\sum_{\delta_i\in F}\int\limits_{
(F\backslash {\bf
A})^2}\theta^{V_1,\psi_U}(l(0,r_2,0,r_4,0,0,0)wg) \psi(\delta_1
r_2+\delta_2 r_4)dr_2dr_4$$ Collapsing summation with integration
this also equals
$$\int\limits_{{\bf
A}^2}\theta^{V_2,\psi_U}(z(m_1,m_2)wg)dm_1dm_2$$ Here $V_2$ is the
unipotent subgroup of $U$ generated by the group $U_1$ and by the
unipotent group consisting of all matrices
$l(r_1,r_2,r_3,r_4,r_5,r_6,r_7)$. We view $\psi_U$ as a character
of $V_2$ by restriction. Performing the same expansion for
$\varphi_\pi$, integral \eqref{global5} equals
\begin{equation}\label{global6}
\int\int\limits_{{\bf
A}^4}\varphi_\pi^{V_2,\psi_U^{-1}}(z(m_3,m_4)wg)
\theta^{V_2,\psi_U}(z(m_1,m_2)wg)f_s(g)dm_idg
\end{equation}
where the $g$ integration is as in \eqref{global5}. Next we expand
$\theta^{V_2,\psi_U}$ along the unipotent subgroup of $U$
generated by all matrices of the form
$t(r_1,r_2,r_3,r_4)=I_{10}+r_1e'_{3,4}+r_2e'_{3,5}+r_3e'_{3,6}+r_4e'_{3,7}$.
It is not hard to check that the only nonzero contribution to the
expansion comes from the constant term. Indeed, all other terms in
the expansion will produce  Fourier coefficients which correspond
to unipotent classes which are greater than $(3^31)$. By
assumption 1) these Fourier coefficients are zero. Let $V_3$
denote the unipotent subgroup of $U$ generated by $V_2$ and all
matrices of the form $t(r_1,r_2,r_3,r_4)$. The above discussion
implies that $\theta^{V_2,\psi_U}=\theta^{V_3,\psi_U}$. Factoring
the integration over the group $Y$, integral \eqref{global6}
equals
\begin{equation}\label{global7}
\int\int\limits_{{\bf A}^4}\int\limits_{(F\backslash {\bf A})^2}
\varphi_\pi^{V_2,\psi_U^{-1}}(t(0,r_2,0,r_4)z(m_3,m_4)wg)
\theta^{V_3,\psi_U}(z(m_1,m_2)wg)f_s(g)dr_idm_jdg
\end{equation}
Here the $g$ variable is integrated over $Z({\bf A})GL_2(F)Y({\bf
A})V({\bf A})\backslash G({\bf A})$.

We now expand the function $\varphi_\pi^{V_2,\psi_U^{-1}}$ in the
above integral along the matrices $t(r_1,0,r_3,0)$ with points in
$F\backslash {\bf A}$. The group $GL_2$ acts on this expansion
with two orbits. The trivial one contributes zero to the integral
because of the cuspidality of $\pi$. Thus \eqref{global7} equals
\begin{equation}\label{global8}
\int\limits_{Z({\bf A})GL_1(F)N'(F)Y({\bf A})V({\bf A})\backslash
G({\bf A})}\int\limits_{{\bf A}^4}
\varphi_\pi^{V_3,{\psi'_U}^{-1}}(z(m_3,m_4)wg)
\theta^{V_3,\psi_U}(z(m_1,m_2)wg)f_s(g)dm_jdg
\end{equation}
Here $N'$ is the maximal unipotent subgroup of $GL_2$ which, when
embedded inside $G$, is a subgroup of $U$. The group $GL_1$
consists of all diagonal matrices inside the group $G$ of the form
$diag(a,a,a,1,1,a,a,1,1,1)$. The character $\psi'_U$, which was
defined before the theorem, is viewed as a character of $V_3$ by
restriction. Next we factor the integration in \eqref{global8}
over the group $N$. The group $GL_1$ acts on this group with two
orbits. The trivial one contributes zero by cuspidality. Thus,
\eqref{global8} equals
\begin{equation}\label{global9}
\int\limits_{Z({\bf A})U_0({\bf A})\backslash G({\bf
A})}\int\limits_{{\bf A}^4} W_\pi(z(m_3,m_4)wg)
\theta^{U,\psi_U}(z(m_1,m_2)wg)f_s(g)dm_jdg
\end{equation}
Here $U_0$ is the group generated by $V$,$Y$ and $N$. Next we
conjugate the matrix $z(m_3,m_4)$ across $w$ and collapse
summation with integration. We also change variables $g\mapsto
w^{-1}g$. Thus  integral \eqref{global9} equals
\begin{equation}\label{global10}
\int\limits_{Z({\bf A})U_1({\bf A})\backslash G({\bf
A})}\int\limits_{{\bf A}^2} W_\pi(g)
\theta^{U,\psi_U}(z(m_1,m_2)wg)f_s(w[53]g)dm_jdg
\end{equation}
where $U_1$ is the subgroup of $U$ defined as follows. Let $U'_0$
be the subgroup of $U_0$ where we omit the one dimensional
unipotent subgroups corresponding to the roots $\alpha_5$ and
$\alpha_3+\alpha_5$. Then $U_1=w^{-1}U'_0w$. Factoring the
integration $U_1({\bf A})\backslash U({\bf A})$, we obtain
integral \eqref{th}. \hfill $\blacksquare$

\section{ A Construction of a Small Representation}

In this section we  construct a representation of the group
$G({\bf A})$ and show that it satisfies the assumptions stated in
Theorem 1. This representation will depend on a choice of a
cuspidal representation of $GL_2({\bf A})$.

Let $\tau$ denote an irreducible cuspidal representation of the
group $GL_2({\bf A})$. We will assume that it has a trivial
central character. Let $\sigma(\tau)$ denote the symmetric square
lift of $\tau$ to $GL_3({\bf A})$. This lift was constructed by
Gelbart and Jacquet in \cite{G-J}. Let $\mu(\tau)=\tau\otimes\tau$
denote the tensor product representation of $GSO_4({\bf A})$.
Denote by $Q$ the standard maximal parabolic subgroup of $G$ whose
Levi part is $GL_3\times GSO_4$. We shall denote by $U(Q)$ its
unipotent radical. Let $E_\tau(g,s)$ denote the Eisenstein series
defined on the group $G({\bf A})$ which is associated to the
induced representation $Ind_{Q({\bf A})}^{G({\bf
A})}(\sigma(\tau)\otimes \mu(\tau))\delta_Q^s$.

From the  Langlands theory the poles of this Eisenstein series are
determined by the poles of the constant terms. Since we induce
from cuspidal data, it follows that we only need to consider the
constant term along $U(Q)$. This is easily computed and hence the
poles of $E_\tau(g,s)$ are determined by
$$\frac{L^S(\sigma(\tau)\times
\mu(\tau),6s-3)L^S(\sigma(\tau),12s-6)}{L^S(\sigma(\tau)\times
\mu(\tau),6s-2)L^S(\sigma(\tau),12s-5)}$$ where $S$ is a finite
set of places, including the archimedean ones, such that outside
of $S$ all data is unramified. The above $L$ function is equal to
$$\frac{L^S(\sigma(\tau)\times
\sigma(\tau),6s-3)L^S(\sigma(\tau),6s-3)
L^S(\sigma(\tau),12s-6)}{L^S(\sigma(\tau)\times
\sigma(\tau),6s-2)L^S(\sigma(\tau),6s-2)L^S(\sigma(\tau),12s-5)}$$
from which we deduce that the Eisenstein series has a simple pole
at $s=2/3$. We denote $\theta_\tau(g)=Res_{s=2/3}E_\tau(g,s)$. In
the rest of this section we will show that this representation
does satisfy the assumptions of Theorem 1.

Since $\sigma(\tau)$ and $\mu(\tau)$ are both generic, the
integral $\theta_\tau^{U,\psi_U}(g)$, which was defined at the
beginning of section 2,  is not zero for some choice of data.
Since this integration corresponds to a unipotent orbit of the
type $(3^31)$, all we need to verify is that $\theta_\tau$ has no
nonzero Fourier coefficient which corresponds to any unipotent
class which is greater than or not related to $(3^31)$.

To do that we use the same method as in \cite{G-R-S2} section two.
We start by studying the unramified local representation
corresponding to $\theta_\tau$. Let $F$ be a local nonarchimedean
field where $\tau$ is unramified. Assume that $\tau$ is a
constituent of the induced representation
$Ind_{B_2}^{GL_2}\chi\delta_{B_2}^{1/2}$, where $\chi$ is an
unramified character of $F^*$ and $B_2$ is the Borel subgroup of
$GL_2$. By the definition of the symmetric square lift, the
unramified local representation corresponding to $\sigma(\tau)$ is
$Ind_{B_3}^{GL_3}\chi_1\delta_{B_3}^{1/2}$. Here $B_3$ is the
Borel subgroup of $GL_3$, and $\chi_1$ is defined as
$\chi_1(diag(a,b,c))=\chi^2(ac^{-1})$. In a similar way we have
$\mu(\tau)=Ind_{B_4}^{GSO_4}\chi_2\delta_{B_4}^{1/2}$. Here $B_4$
is the Borel subgroup of $GSO_4$ and
$\chi_2(diag(abr,ar,a^{-1},a^{-1}b^{-1}))=\chi^2(ab)\chi(r)$.
Here the factor $\chi(r)$ occurs because we assume that
$\mu(\tau)$ has a trivial central character.

From all this we deduce\\
{\bf Lemma 2:} {\sl Let $\theta'_\tau$ denote the unramified
constituent of $\theta_\tau$ at a nonarchimedean place. Then
$\theta'_\tau$ is a sub-quotient of
$Ind_{P}^{G}\chi_3\delta_Q^{1/2}$. Here, for all $g\in GL_3$ and
$h\in GSO_4$ we define
$\chi_3((g,h))=\chi^2(detg)\chi^3(\lambda(h))$
where $\lambda(h)$ denotes the similitude factor of the matrix $h$.}\\
{\sl Proof:} Let $B$ denote the Borel subgroup of $G$. The
unramified representation $\theta'_\tau$ is a constituent of the
induced representation $Ind_B^G\chi_4\delta_B^{1/2}$ where
$\chi_4(t)=\chi^2(a_1a_3^{-1}a_4)\chi(r)$ and
$t=diag(ra_1,ra_2,ra_3,ra_4,r,1,a_4^{-1},a_3^{-1},a_2^{-1},a_1^{-1})$.
Let $w_0$ denote the Weyl element of $G$ which has a one at the
entries $$(1,1);(2,4);(3,8);(4,2);(5,6);(6,5);
(7,9);(8,3);(9,7);(10,10)$$ and zero elsewhere. Then
$Ind_B^G\chi_4\delta_B^{1/2}$ is isomorphic to
$Ind_B^G\chi_4^{w_0}\delta_B^{1/2}$ where
$\chi_4^{w_0}(t)=\chi_4(w_0^{-1}tw_0)$. This can be written as
$Ind_B^G(\chi_5\delta_{B_3}^{1/2}\chi_6\delta_{B_4}^{1/2})\delta_Q^{1/2}$
where $\chi_5(t)=\chi^2(a_1a_2a_3)$ and $\chi_6(t)=\chi^3(r)$.
From this it follows that $\theta'_\tau$ is a sub-quotient of
$Ind_{P}^{G}\chi_3\delta_Q^{1/2}$. \hfill $\blacksquare$

We return to the global situation. We need to prove that
$\theta_\tau$ has no nonzero Fourier coefficient which corresponds
to any unipotent class which is greater than or not related to
$(3^31)$. As mentioned in section two, in \cite{G-R-S1} it is
explained how to associate with a unipotent class of a group $G$ a
set of Fourier coefficients. It is clear that $\theta_\tau$ is not
generic. Hence it has no nonzero Fourier coefficient with respect
to the unipotent class $(91)$. Arguing as in \cite{G-R-S1} lemma
2.6 we deduce that if $\theta_\tau$ has a nonzero Fourier
coefficient corresponding to the unipotent class $(73)$ then it
has a nonzero Fourier coefficient which corresponds to the
unipotent class $(71^3)$. Similarly, if $\theta_\tau$ has a
nonzero Fourier coefficient which corresponds to the unipotent
class $(5r_1\ldots r_k)$ then it has a nonzero Fourier coefficient
which corresponds to the unipotent class $(51^5)$. Arguing as in
\cite{G-R-S2} lemma 3, we deduce that if $\theta_\tau$ has no
nonzero Fourier coefficients which corresponds to the unipotent
class $(51^5)$ then it has no nonzero Fourier coefficient
corresponding to the unipotent class $(71^3)$. We should remark
that these last statements are true for any automorphic
representation of $G$ and not only for $\theta_\tau$. Hence it is
enough to show that $\theta_\tau$ has no nonzero Fourier
coefficients which correspond to the unipotent classes $(51^5)$
and to $(4^21^2)$. We now describe the  families of Fourier
coefficients which correspond to these two unipotent classes.

Let $V$ denote the unipotent radical subgroup of the standard
parabolic subgroup of $G$ whose Levi part is $GL_1^2\times GSO_6$.
We view $V$ as a subgroup of the maximal unipotent subgroup $U$ of
$G$ which consists of upper triangular matrices. Let $a\in F^*$.
We define a character $\psi_{V,a}$ of the group $V$ as follows.
For $v=(v_{i,j})\in V$, set
$\psi_{V,a}(v)=\psi(v_{1,2}+v_{2,5}+av_{2,6})$. Thus the Fourier
coefficient of the representation $\theta_\tau$ which corresponds
to the unipotent class $(51^5)$ is given by
$$\int\limits_{V(F)\backslash V({\bf
A})}\theta_\tau(vg)\psi_{V,a}(v)dv$$ Next we describe the Fourier
coefficient which corresponds to the unipotent class $(4^21^2)$.
To do that, let $R$ denote the subgroup of $U$ defined as the
group generated by all one dimensional unipotent subgroup
associated with the positive roots of $G$ where we omit the roots
$\alpha_1,\alpha_3,\alpha_4$ and $\alpha_5$. Thus, the dimension
of $R$ is 16. For $r=(r_{i,j})\in R$ we define
$\psi_R(r)=\psi(r_{1,3}+r_{2,4}+r_{3,7})$. The Fourier coefficient
which corresponds to the unipotent class $(4^21^2)$ is given by
$$\int\limits_{R(F)\backslash R({\bf
A})}\theta_\tau(rg)\psi_{R}(r)dr$$ To prove our result we will
show that the local unramified component $\theta'_\tau$ of
$\theta_\tau$ does not support a local functional of the above type.
More precisely, we prove\\
{\bf Lemma 3:} {\sl Let ${\cal O}$ denote one of the unipotent
orbits $(51^5)$ or $(4^21^2)$. Then $\theta^{'}_\tau$ has no
nonzero linear functional $l_{{\cal O},a}$ which satisfies
$l_{{\cal O},a}(\rho(v)x)=\psi_{V,a}^{-1}(v)l_{{\cal O},a}(x)$ and
no nonzero linear functional $l_{{\cal O}}$ which satisfies
$l_{{\cal O}}(\rho(r)x)=\psi_{R}^{-1}(r)l_{{\cal O}}(x)$. This is
for all $v\in V$, $r\in R$ and $x$ is a vector in the space of
$\theta^{'}_\tau$. Here we denoted by $\rho$ the action of the
representation $\theta^{'}_\tau$.}\\
{\sl Proof:} Arguing as in \cite{G-R-S2} section 2 lemma 2 we have
to show that the representation $Ind_{Q}^{G}\chi_3\delta_Q^{1/2}$
does not support any of these functionals. From the Bruhat theory
this reduces to the problem of showing that
$Ind_{Q}^{G}\chi_3\delta_Q^{1/2}$ has no admissible double coset
in the space $Q\backslash G/L$ where $L$ is either $V$ or $R$. By
that we mean that for any $g\in Q\backslash G/L$ there is $l\in L$
such that $glg^{-1}\in Q$ and that $\psi_L(l)\ne 1$. Here $\psi_L$
is either $\psi_{V,a}$ or $\psi_R$. From the Bruhat decomposition
we obtain that each element $g\in Q\backslash G/L$ can be written
as $g=wu$ where $w$ is a Weyl element of $G$ and $u\in U$. As  in
\cite{G-R-S2} section 2 lemma 2 we deduce that $wu$ is not
admissible if and only if $w$ is not admissible. Let $U(Q)^-$
denote the transpose group of $U(Q)$. Thus $U(Q)^-$ consists of
lower unipotent matrices. Let $x_1(r_1)=I_{10}+r_1e'_{1,2},
x_2(r_2)=I_{10}+r_2e'_{2,5}$ and $x_3(r_3)=I_{10}+r_3e'_{2,6}$.
These are precisely the three one dimensional unipotent groups on
which the character $\psi_{V,a}$ is not trivial. To prove that $w$
is not admissible it is enough to show that for some $i$, we have
$wx_i(r_i)w^{-1}\in Q$. Assume not. This means that
$wx_i(r_i)w^{-1}\in U(Q)^-$ for all $i$. Since $w$ is in $G$ we
have that if $w_{i,j}=1$ then $w_{11-i,11-j}=1$. This means that
if $wx_2(r_2)w^{-1}=I_{10}+r_2e'_{i,j}$ then
$wx_3(r_3)w^{-1}=I_{10}+r_3e'_{11-i,j}$. Hence if
$wx_2(r_2)w^{-1}\in [U(Q)^-,U(Q)^-]$ then $wx_3(r_3)w^{-1}\in Q$.
Similarly, if $wx_3(r_3)w^{-1}\in [U(Q)^-,U(Q)^-]$ then
$wx_2(r_2)w^{-1}\in Q$.  From this it follows that both
$wx_2(r_2)w^{-1},wx_3(r_3)w^{-1}\in U(Q)^-/[U(Q)^-,U(Q)^-]$. Since
$x_1(r_1)$ does not commute with both $x_2(r_2)$ and $x_3(r_3)$,
it follows the same for $wx_i(r_i)w^{-1}$. Hence, also
$wx_1(r_1)w^{-1}\in U(Q)^-/[U(Q)^-,U(Q)^-]$. However, there is no
one dimensional unipotent subgroup of $G$ which corresponds to a
root in $G$, which has the property that it is in
$U(Q)^-/[U(Q)^-,U(Q)^-]$ and which will commute with both
$wx_2(r_2)w^{-1}$ and $wx_3(r_3)w^{-1}$. Thus every $w$ is not
admissible and the lemma is proved for the functional $l_{{\cal
O},a}$. A similar proof holds for the group $R$ and the functional
$l_{{\cal O}}$. \hfill $\blacksquare$

\section{ The Unramified Computation}\label{localsection}

In this section we carry out the unramified computation of the
local integral which corresponds to integral \eqref{th}. Let $F$
denote a nonarchimedean field, and assume that all groups are
defined over $F$. Let $\pi$ denote an unramified irreducible
generic representation of the group $G$. We assume that it has a
trivial central character. We denote by $I(s)$ the local induced
representation $Ind_P^G\delta_P^s$. Let $\tau$ denote an
unramified irreducible representation of $GL_2$ and assume that it
has a trivial central character. As in section 3 lemma 2 we denote
by $\theta'_\tau$ the local unramified constituent of the global
representation $\theta_\tau$. The integral we consider is
\begin{equation}\label{local1}
\int\limits_{ZU\backslash G}\int\limits_{F^2}
W_\pi(g)\theta'^{U,\psi_U}_\tau(g)
f_s(w[53]x_{\alpha_3}(r_1)x_{\alpha_3+\alpha_5}(r_2)g)\psi(r_1)dr_1dr_2dg
\end{equation}
Assume that all functions are the $K$ fixed vectors in their space
of representation, where $K$ is the maximal compact subgroup of
$G$. Let $Spin_{10}$ denote the 16 dimensional irreducible Spin
representation of $GSpin_{10}({\bf C})$, the $L$ group of $G$. The
corresponding $L$ function was defined in \cite{G1}. We denote by
$L(Spin_{10}\times St,\pi\times \tau,s)$ the 32 dimensional $L$
function which consists of the tensor product of these two
representations. In this section we prove\\
{\bf Proposition 4:} {\sl For all unramified data, and for $Re(s)$
large, integral \eqref{local1} equals
\begin{equation}
\frac{L(Spin_{10}\times St,\pi\times
\tau,2s-1/2)}{\zeta(8s)\zeta(8s-2)}
\end{equation}}
{\sl Proof:} Let $T$ denote the maximal torus of $G$. Using the
Iwasawa decomposition, integral \eqref{local1} equals
$$\int\limits_{Z\backslash
T}W_\pi(t)\theta'^{U,\psi_U}_\tau(t)\delta_B^{-1}(t)
\int\limits_{F^2}f_s(w[53]x_{\alpha_3}(r_1)x_{\alpha_3+\alpha_5}
(r_2)t)\psi(r_1)dr_1dr_2dt$$ where $B$ is the Borel subgroup of
$G$. We parameterize an element $t$ in $Z\backslash T$ as
$t=diag(ab_1,ab_2,ab_3,ab_4,a,1,b_4^{-1},b_3^{-1},b_2^{-1},b_1^{-1})$.
For this parameterization
$\delta_B^{-1}(t)=|a|^{-10}|b_1|^{-8}|b_2|^{-6}|b_3|^{-4}|b_4|^{-2}$.

Next we compute the inner integration along the two dimensional
unipotent subgroup. This is done exactly as in \cite{G-H} section
3 right after identity (11). We thus obtain
$$\int\limits_{F^2}f_s(w[53]x_{\alpha_3}(r_1)x_{\alpha_3+\alpha_5}
(r_2)t)\psi(r_1)dr_1dr_2=$$
$$\frac{\zeta(8s-2)}{\zeta(8s)}(1-|b_3b_4^{-1}|^{8s-2}q^{-8s+2})
|a|^{2s+1}|b_1b_2|^{4s}|b_3|^{-4s+2}|b_4|^{4s-1}$$ Here we used
the fact that $\delta_P(t)=|a|^{10s}|b_1b_2b_3b_4|^{4s}$. Also, we
denote $q=|p|$ where $p$ is a generator of the maximal ideal of
the ring of integers in the field $F$.

Recall, from section 3, that $\theta'^{U,\psi_U}_\tau$ is the
local component of a residue of an Eisenstein series at the point
$2/3$. Hence we have
$$\theta'^{U,\psi_U}_\tau(t)=\delta_Q^{1/3}(t)W_\tau\begin{pmatrix}
b_1&&\\&b_2&\\&&b_3 \end{pmatrix}W_{\tau\times\tau}\begin{pmatrix}
ab_4&&&\\&a&&\\&&1&\\&&&b_4^{-1} \end{pmatrix}$$ Here $W_\tau$ is
the unramified Whittaker function corresponding to the induced
representation $Ind_{B_3}^{GL_3}\chi_1\delta_{B_3}^{1/2}$ of
$GL_3$ which was defined in section 3. Similarly, the function
$W_{\tau\times\tau}$ is the unramified Whittaker function
corresponding to the induced representation
$Ind_{B_4}^{GSO_4}\chi_2\delta_{B_4}^{1/2}$. Denote
$K_\pi(t)=W_\pi(t)\delta_B^{-1/2}(t)$ and also
$$K_\tau\begin{pmatrix}
b_1&&\\&b_2&\\&&b_3 \end{pmatrix}K_{\tau\times\tau}\begin{pmatrix}
ab_4&&&\\&a&&\\&&1&\\&&&b_4^{-1}
\end{pmatrix}=$$ $$W_\tau\begin{pmatrix} b_1&&\\&b_2&\\&&b_3
\end{pmatrix}W_{\tau\times\tau}\begin{pmatrix}
ab_4&&&\\&a&&\\&&1&\\&&&b_4^{-1}
\end{pmatrix}|b_1^{-1}b_3||a|^{-1/2}|b_4|^{-1}$$ Using the fact
that $\delta_Q^{1/3}(t)=|a|^3|b_1b_2b_3|^2$ we obtain that
integral \eqref{local1} equals
$$\frac{\zeta(8s-2)}{\zeta(8s)}\int\limits_{Z\backslash T}K_\pi(t)
K_\tau\begin{pmatrix} b_1&&\\&b_2&\\&&b_3
\end{pmatrix}K_{\tau\times\tau}\begin{pmatrix}
ab_4&&&\\&a&&\\&&1&\\&&&b_4^{-1}
\end{pmatrix}\times$$
$$(1-|b_3b_4^{-1}|^{8s-2}q^{-8s+2})|a|^{2s-1/2}|b_1b_2b_3^{-1}b_4|^{4s-1}dt$$
We consider the following change of variables. Set $a\mapsto
t_1t_5^{-1}, b_1\mapsto t_2t_3t_4t_5, b_2\mapsto t_3t_4t_5,
b_3\mapsto t_4t_5$ and $b_4\mapsto t_5$. With this change of
variables the torus $Z\backslash T$ is parameterized as
$t=diag(t_1t_2t_3t_4t_5,t_1t_3t_4t_5,t_1t_4t_5,t_1t_5,t_1,t_5,1,t_4^{-1},
t_3^{-1}t_4^{-1},t_2^{-1}t_3^{-1}t_4^{-1})$. Thus the above
integral equals
$$\frac{\zeta(8s-2)}{\zeta(8s)}\int\limits_{Z\backslash T}K_\pi(t)
K_\tau\begin{pmatrix} t_2t_3&&\\&t_2&\\&&1
\end{pmatrix}K_{\tau\times\tau}\begin{pmatrix}
t_1t_5&&&\\&t_1&&\\&&t_5&\\&&&1
\end{pmatrix}\times$$
$$(1-|t_4|^{8s-2}q^{-8s+2})|t_1|^{2s-1/2}
|t_2t_4|^{4s-1}|t_5|^{6s-3/2}|t_3|^{8s-2}dt$$ For $1\le i\le 5$
write $t_i=p^{n_i}$. Also, we set $y=q^{-2s+1/2}$. Denote by
$(m_1,m_2,m_3,m_4,m_5)$ the trace of the irreducible
representation $\sum{m_i}\varpi_i$ evaluated at the semi-simple
conjugacy class of $GSpin_{10}({\bf C})$ associated with $\pi$.
Here $\varpi_i$ is the i-th fundamental representation of
$GSpin_{10}({\bf C})$. It follows from the Casselman-Shalika
formula \cite{C-S} that $K_\pi(t)=(n_2,n_3,n_4,n_5,n_1)$.
Similarly, we denote by $(m)$ the trace of the irreducible
representation $m\varpi$ evaluated at the semi-simple conjugacy
class of $GL_{2}({\bf C})$ associated with $\tau$. Using again the
Casselman-Shalika formula we obtain
$$K_\tau\begin{pmatrix} t_2t_3&&\\&t_2&\\&&1
\end{pmatrix}K_{\tau\times\tau}\begin{pmatrix}
t_1t_5&&&\\&t_1&&\\&&t_5&\\&&&1
\end{pmatrix}=l(n_2,n_3)\otimes (n_1)\otimes (n_5)$$ where
$l(n_2,n_3)$ denotes the restriction of the irreducible
representation with highest weight $[n_2,n_3]$ of $GL_3({\bf C})$
to the group $SO_3({\bf C})$.

Thus, to prove the proposition we are reduced to prove the
identity
$$\sum_{n_i}^{\infty}(n_2,n_3,n_4,n_5,n_1)k(n_1,n_2,n_3,n_5)
y^{n_1+2n_2+4n_3+2n_4+3n_5}(1-y^{4(n_4+1)})=$$
$$(1-y^4)^2L(Spin_{10}\times St,\pi\times \tau,2s-1/2)$$ where
$k(n_1,n_2,n_3,n_5)=l(n_2,n_3)\otimes (n_1)\otimes (n_5)$. For
$\chi$ as in section 3, we have
$$L(Spin_{10}\times St,\pi\times \tau,2s-1/2)
=L(Spin_{10},\pi\otimes \chi,2s-1/2)
L(Spin_{10},\pi\otimes \chi^{-1},2s-1/2).
$$
By the Poincar\'e identity, and Brion's decomposition of the
symmetric algebra of $Spin_{10}$  ( see \cite{Br}) we have
$$L(Spin_{10},\pi\otimes \chi,2s-1/2)=\sum_{m,\ell=0}^{\infty}
(m,0,0,0,\ell)
\chi^{2m+\ell}y^{2m+\ell}.$$
Let
$$P(\chi,y)=(1-\chi^8y^8)+(\chi^6y^6-\chi^2y^2)(1,0,0,0,0)
\chi^2y^2+
\chi^3y^3(0,0,0,1,0)-\chi^5y^5(0,0,0,0,1).$$ Then, using again the
Casselman-Shalika formula, we obtain
$$P(\chi,y)
\sum_{m,\ell=0}^{\infty}(m,0,0,0,\ell)
\chi^{2m+\ell}y^{2m+\ell}=\sum_{\ell=0}^{\infty}(0,0,0,0,\ell)
\chi^{\ell}y^{\ell}.$$

So, our identity reduces to
$$P(\chi,y)P(\chi^{-1},y)
\sum_{n_i}^{\infty}(n_2,n_3,n_4,n_5,n_1)k(n_1,n_2,n_3,n_5)
y^{n_1+2n_2+4n_3+2n_4+3n_5}(1-y^{4(n_4+1)})$$
$$=(1-y^4)^2\sum_{m,\ell=0}^{\infty}
(0,0,0,0,\ell)\otimes (0,0,0,0,m)\chi^{m-\ell}y^{m+\ell}.$$ Here
and henceforth, by abuse of notation, we write
$(n_2,n_3,n_4,n_5,n_1)$ for the representation itself, as well as
the value of the trace at the conjugacy class associated to $\pi.$

{\bf Lemma 5:} {\sl We have
$$
(0,0,0,0,\ell)\otimes (0,0,0,0,m)=\sum_{a+b\leq min(m,\ell)}
(a,0,b,0,m+\ell-2a-2b).$$} 
{\sl Proof:}  This result is due to Okada \cite{O}
see also \cite{Kr}. \hfill $\blacksquare$

What remains is to check that
$$P(\chi,y)P(\chi^{-1},y)
\sum_{n_i}^{\infty}(n_2,n_3,n_4,n_5,n_1)k(n_1,n_2,n_3,n_5)
y^{n_1+2n_2+4n_3+2n_4+3n_5}(1-y^{4(n_4+1)})$$
\begin{equation}\label{pid}
=(1-y^4)^2\sum_{m,\ell=0}^{\infty}
\left(\sum_{n_2+n_4\leq min(m,\ell)}
(n_2,0,n_4,0,m+\ell-2n_2-2n_4)\right)\chi^{m-\ell}y^{m+\ell}
\end{equation}
$$=(1-y^4)^2\sum_{n_2,n_4,n_1=0}^{\infty}(n_2,0,n_4,0,n_1)
y^{2n_2+2n_4+n_1}\chi^{-n_1}\frac{1-\chi^{2(n_1+1)}}{1-\chi^2}.$$
We check \eqref{pid} in two stages.  First,
we show that
\begin{equation}\label{pchi}
P(\chi,y)
\sum_{n_i}^{\infty}(n_2,n_3,n_4,n_5,n_1)k(n_1,n_2,n_3,n_5)
y^{n_1+2n_2+4n_3+2n_4+3n_5}(1-y^{4(n_4+1)})\end{equation}
$$=(1-y^4)^2\sum_{n_i}^{\infty}\frac{1-\chi^{2(n_1+1)}}{1-\chi^2}
\frac{1-\chi^{2(n_2+1)}}{1-\chi^2}\chi^{-(n_1+2n_2+2n_3+n_5)}
y^{n_1+2n_2+4n_3+2n_4+3n_5}.$$ Then, we check that
$P(\chi^{-1},y)$ times this sum is the desired sum. First, using
the Casselman-Shalika formula, we note that $k(n_1,n_2,n_3,n_5)$
is equal to
$$\frac{(1-\chi^{2(n_2+1)})(1-\chi^{2(n_3+1)})
(1-\chi^{2(n_2+n_3+2)})(1-\chi^{2(n_1+1)})(1-\chi^{2(n_5+1)})}
{(1-\chi^{2})^4(1-\chi^4)}\chi^{-(n_1+2n_2+2n_3+n_5)}.$$ To check
\eqref{pchi}, we shall verify that the coefficient of
$(n_2,n_3,n_4,n_5,n_1)$ is the same on both sides.  Let
$$g_{[n_i]}(\chi,y)=\frac{(1-\chi^{2(n_2+1)})(1-\chi^{2(n_3+1)})
(1-\chi^{2(n_2+n_3+2)})(1-\chi^{2(n_1+1)})(1-\chi^{2(n_5+1)})}
{(1-\chi^{2})^4(1-\chi^4)}$$
$$\chi^{-(n_1+2n_2+n_3+n_5)}y^{n_1+2n_2+4n_3+2n_4+3n_5}
(1-y^{4(n_4+1)}).$$
Let $\Gamma_1$ denote the set of weights of $(1,0,0,0,0)$,
$\Gamma_4$ those of $(0,0,0,1,0)$, and $\Gamma_5$ those of
$(0,0,0,0,1)$. Then, we claim that
$$P(\chi,y)
\sum_{n_i}^{\infty}(n_2,n_3,n_4,n_5,n_1)g_{[n_i]}(\chi,y)=
\sum_{n_i}^{\infty}(n_2,n_3,n_4,n_5,n_1)G_{(n_i)}(\chi,y)$$
where
\begin{equation}
G_{(n_i)}(\chi,y)=
(1-\chi^8y^8)g_{[n_i]}(\chi,y)+(\chi^6y^6-\chi^2y^2)
\sum\limits_{w \in
\Gamma_1}g_{[n_i-w_i]}(\chi,y)+\label{gni}\end{equation}
$$\chi^3y^3\sum\limits_{w \in \Gamma_4}g_{[n_i-w_i]}(\chi,y)
-\chi^5y^5\sum\limits_{w \in \Gamma_5}g_{[n_i-w_i]}(\chi,y)).$$
Indeed, Brauer proved in general that if $\chi_{\lambda}$ denotes
the character
of the irreducible representation of highest weight $\lambda,$ then
$$\chi_{\lambda}\chi_{\mu}=\sum_{\nu}s_\nu
\chi_{|\lambda+\nu+\rho|-\rho},$$
where
the sum is over the weights of the representation with highest
weight $\mu,$ with multiplicity,
$|\tau|$ denotes the dominant member of the Weyl orbit of
$\tau$,
and $s_\nu$ is equal to $(-1)^w$ if
$w(\lambda+\nu+\rho)=|\lambda+\nu+\rho|$
for a unique $w$, and zero if $\lambda+\nu+\rho$ has a stabilizer
 in the
Weyl group.  See \cite{Bu}
p. 171 and exercises.  In our case, no element of
$\Gamma_1\cup\Gamma_4\cup\Gamma_5$ has an entry less than $-1$, so
when $\lambda+\nu+\rho$ is not dominant, it has a stabilizer in
the Weyl group.  This proves \eqref{gni} when no $n_i$ is zero. If
any of the $n_i$ is zero, there are extra terms on the right hand
side corresponding to $n_i-w_i=-1.$  But referring back to the
formula for $g_{[n_i]}$, we see that all these terms vanish
anyway.

Hence, to check \eqref{pchi}, we now need to check that
$G_{(n_i)}$ is equal to the coefficient of $(n_2,n_3,n_4,n_5,n_1)$
on the right side of \eqref{pchi}. Replacing $\chi^{2(n_i+1)}$ by $X_i$
and $y^{4n_4}$ by $Y_4$ we obtain an identity of polynomials in seven 
variables which is straightforward to verify by computer.

Now, we need to check that
$$P(\chi^{-1},y)\sum_{n_i}^{\infty}
\frac{1-\chi^{2(n_1+1)}}{1-\chi^2}
\frac{1-\chi^{2(n_2+1)}}{1-\chi^2}\chi^{-(n_1+2n_2+2n_3+n_5)}
y^{n_1+2n_2+4n_3+2n_4+3n_5}$$
$$=\sum_{n_2,n_4,n_1=0}^{\infty}
(n_2,0,n_4,0,n_1)
y^{2n_2+2n_4+n_1}\chi^{-n_1}\frac{1-\chi^{2(n_1+1)}}{1-\chi^2}.$$
Let $$h_{[n_i]}(\chi,y)=\frac{1-\chi^{2(n_1+1)}}{1-\chi^2}
\frac{1-\chi^{2(n_2+1)}}{1-\chi^2}\chi^{-(n_1+2n_2+2n_3+n_5)}
y^{n_1+2n_2+4n_3+2n_4+3n_5}$$ when all $n_i$ are nonnegative. Note
that this expression is equal to zero if $n_1$ or $n_2$ is equal
to $-1$. We extend the definition of $h_{[n_i]}$ to be zero if
$\min(n_3,n_4,n_5)=-1.$ Then the left hand side is equal to
$$\sum_{n_i}^{\infty}(n_2,n_3,n_4,n_5,n_1)H_{(n_i)}(\chi,y),$$
where
$$H_{(n_i)}(\chi,y)=
(1-\chi^{-8}y^8)h_{[n_i]}(\chi,y)+(\chi^{-6}y^6-\chi^{-2}y^2)
\sum\limits_{w \in \Gamma_1}h_{[n_i-w_i]}(\chi,y)+$$
$$\chi^{-3}y^3\sum\limits_{w \in \Gamma_4}h_{[n_i-w_i]}(\chi,y)
-\chi^{-5}y^5\sum\limits_{w \in \Gamma_5}h_{[n_i-w_i]}(\chi,y)).$$
We claim that this is equal to
$$
y^{2n_2+2n_4+n_1}\chi^{-n_1}\frac{1-\chi^{2(n_1+1)}}{1-\chi^2},$$
when $n_3=n_5=0$, and zero otherwise.
To show this, we partition each $\Gamma_i$ into eight sets,
according to
whether the $w_3,w_4,$ and $w_5$ are positive or nonpositive.  This
determines in which $H_{(n_i)}$ the corresponding term is zero by
convention.

We replace $\chi^{2(n_1+1)}$ and $\chi^{2(n_2+1)}$ by $X_1$ and
$X_2$ throughout. A factor of
$$\chi^{-(n_1+2n_2+2n_3+n_5)}
y^{n_1+2n_2+4n_3+2n_4+3n_5}(1-\chi^2)^{-2}$$ may be factored out
of the sum and moved to the other side,  which
 then
becomes
$$R(\chi,X_1,X_2)=
\frac{X_2}{\chi^2}(1-X_1)(1-\chi^2),$$
when $n_3=n_5=0.$
On the other side of the equation we have a polynomial obtained by
summing
$$(1-X_1\chi^{-2w_1})(1-X_2\chi^{-2w_2})\chi^{w_1+2w_2+2w_3+w_5}
y^{-(w_1+2w_2+4w_3+2w_4+3w_5)}$$ over the weights, with the correct
coefficients from $P$.  We break this into eight pieces.

Let $Q_{nnn}$ denote the sum over weights where $w_3,w_4,w_5$ are
all
nonpositive, $Q_{ppp}$ the sum where they are all positive,
$Q_{npn}$
the sum over weights such that $w_3$ and $w_5$ are nonpositive,
and $w_4$ is positive, and so on.  Since
$(1-\chi^8y^8)h_{[n_i]}(\chi,y)$ is there for every $[n_i]$, it is
 included
in $Q_{nnn}.$  The answers are:
$$Q_{ppp}=Q_{ppn}=Q_{npp}=Q_{npn}=0$$
$$Q_{nnn}=Q_{pnp}=R$$
$$Q_{nnp}=Q_{pnn}=-R.$$
From here it is an easy
application of the inclusion-exclusion principle. (It is easy to
check by hand that $Q_{ppp}=Q_{ppn}=Q_{npp}=0,$ since in all three
cases the set of roots satisfying the desired condition is
empty.)

\section{\bf A construction of a certain lifting}

In this section we shall use a certain small representation,
defined on the group $SO_{22}({\bf A})$, in order to construct the
endoscopic lifting from $PGL_2\times Sp_6$ to the group $SO_{10}$.
Then, in the next section we will use this construction to prove
our main result.

This construction is a special case of a more general set of
constructions which are defined on classical groups. It is a part
of a work in progress of the first named author and was partly
announced in \cite{G2}. To make things clear, we shall now sketch
this construction and just state the main results. The details of
the proofs  will appear in \cite{G3}.

Let $\tau$ denote a cuspidal representation of $PGL_2({\bf A})$.
In section 3 we constructed a small representation $\theta_\tau$,
defined on the group $GSO_{10}({\bf A})$. Clearly the construction
there is also valid if we consider the group $SO_{10}$ instead.
Recall that $\sigma(\tau)$ is the cuspidal representation of
$GL_3({\bf A})$ obtained by the symmetric square lift of $\tau$,
as constructed in \cite{G-J}. Let $E_\tau(g,s)$ denote the
Eisenstein series defined on $GL_6({\bf A})$ associated with the
induced representation $Ind_{L({\bf A})}^{GL_6({\bf
A})}(\sigma(\tau)\otimes\sigma(\tau))\delta_L^s$. Here $L$ is the
maximal parabolic subgroup of $GL_6$ whose Levi part is
$GL_3\times GL_3$. It is well known, see for example \cite{J},
that this representation has a unique simple pole and its residue,
which we shall denote by $E_\tau$, is the well known  Speh
representation. Let $Q$ denote the maximal parabolic subgroup of
the split orthogonal group $SO_{22}$, whose Levi part is
$GL_6\times SO_{10}$. Let $\widetilde{E}_\tau(m,s)$ denote the
Eisenstein series defined on $SO_{22}({\bf A})$ associated with
the induced representation $Ind_{Q({\bf A})}^{SO_{22}({\bf
A})}(E_\tau\otimes\theta_\tau)\delta_Q^s$.

The following two facts about this Eisenstein series are proved in
\cite{G3}.\\
{\bf a)} In the domain $Re(s)>1/2$, the Eisenstein series
$\widetilde{E}_\tau(m,s)$ has a unique simple pole. We shall
denote the residue representation which is obtained by
$\widetilde{\theta}_\tau$. By abuse of notations we shall denote
by $\widetilde{\theta}_\tau(m)$ a vector in this representation
when realized in the space of automorphic
forms.\\
{\bf b)} Using the notations introduced in the beginning of
section 2, we have ${\cal
O}_{SO_{22}}(\widetilde{\theta}_\tau)=(3^71)$. This means that
$\widetilde{\theta}_\tau$ has a nonzero Fourier coefficient which
corresponds to the unipotent orbit $(3^71)$. Also, for any
unipotent orbit ${\cal O}$ which is greater than or not related to
$(3^71)$, the representation $\widetilde{\theta}_\tau$ has no
nonzero Fourier coefficient corresponding to the orbit ${\cal O}$.

At this point we are ready to introduce the global lifting we need
for our result. Let $\epsilon_1$ denote  a cuspidal representation
defined on the group $Sp_6({\bf A})$. Let
$\widetilde{\theta}_\phi$ denote the theta representation defined
on the group $\widetilde{Sp}_{60}({\bf A})$. Here $\widetilde{Sp}$
denotes the double cover of the symplectic group and $\phi$
denotes a Schwartz function defined on ${\bf A}^{30}$. Let $U$
denote the unipotent radical subgroup of $Q$ and let $H_{61}$
denote the Heisenberg group with $61$ variables. We shall write an
element $h\in H_{61}$ as $h=(X,Y,z)$ where $X,Y\in Mat_{3\times
10}$ and $z\in Mat_{1\times 1}$. We define a homomorphism
$l:U\mapsto H_{31}$ as follows. In terms of matrices we can
identify $U$ with the group of all matrices
\begin{equation}\label{mat}
\left \{ u=\begin{pmatrix} I_3&&Y&R_1&R_2\\&I_3&X&R_3&R_1^*\\
&&I_{10}&X^*&Y^*\\ &&&I_3&\\ &&&&I_3 \end{pmatrix} \right \}
\end{equation}
Here $X,Y\in Mat_{3\times 10}$, $R_1\in Mat_{3\times 3}$ and
$R_2,R_3\in Mat^0_{3\times 3}=\{ R\in Mat_{3\times 3}\ :\
R^tJ_3+J_3R=0\}$. The matrix $J_3$ was defined in the beginning of
section 2. For all $u\in U$ as above, if $R_1=(r_{i,j})$ we define
$l(u)=(X,Y,r_{1,1}+r_{2,2}+r_{3,3})$. This map is clearly onto.

We now consider
\begin{equation}\label{lift1}
f_1(g)=\int\limits_{Sp_6(F)\backslash Sp_6({\bf A})}\int
\limits_{U(F)\backslash U({\bf A})}\varphi_{\epsilon_1}(h)
\widetilde{\theta}_\phi(l(u)(h,g))\widetilde{\theta}_\tau(u(h,g))dudh
\end{equation}
Here $g\in SO_{10}({\bf A})$, and the embedding of the groups
$Sp_6$ and $SO_{10}$ in $Sp_{60}$ and $SO_{22}$ are as follows.
First, inside $Sp_{60}$ we embed these two groups  via the tensor
product representation. In $SO_{22}$, we embed them in the obvious
way inside the Levi part of $Q$ which is $GL_6\times SO_{10}$.

Let $\pi_1$ denote the automorphic representation on $SO_{10}({\bf
A})$ defined by all right translations of the functions $f_1(g)$
defined in \eqref{lift1}. We summarize some of the properties of
this representation. The proofs are given in \cite{G3}.\\
{\bf 1)} The representation $\pi_1$ is a cuspidal representation
provided a certain integral is zero. This is the tower property,
and the certain integral we refer to is the lifting in the
previous stage. Since we will not need this point we shall not
discuss it any further.\\
{\bf 2)} The representation $\pi_1$ is generic if and only if
$\epsilon_1$ is a generic cuspidal representation of $Sp_6({\bf
A})$.

It is  important to us that we can interchange the representations
$\pi_1$ and $\epsilon_1$ in integral \eqref{lift1}. More
precisely, let now $\pi$ denote a cuspidal representation defined
on $SO_{10}({\bf A})$. Then we can form the space of automorphic
functions defined by
\begin{equation}\label{lift2}
f(h)=\int\limits_{SO_{10}(F)\backslash SO_{10}({\bf A})}\int
\limits_{U(F)\backslash U({\bf A})}\varphi_\pi(g)
\widetilde{\theta}_\phi(l(u)(h,g))\widetilde{\theta}_\tau(u(h,g))dudg
\end{equation}
Let $\epsilon$ denote the automorphic representation of $Sp_6({\bf
A})$ defined by right translations of the above functions $f(h)$.
The last fact that we need, and which will be proved in \cite{G3},
is the relations at the unramified components of each of these
representations. More precisely, let $\pi,\epsilon$ and $\tau$ be
as above. Suppose that the integral
\begin{equation}\label{howe}
\int\limits_{Sp_6(F)\backslash Sp_6({\bf
A})}\int\limits_{SO_{10}(F)\backslash SO_{10}({\bf A})}\int
\limits_{U(F)\backslash U({\bf A})}f_\epsilon(h)\varphi_\pi(g)
\widetilde{\theta}_\phi(l(u)(h,g))\widetilde{\theta}_\tau(u(h,g))dudgdh
\end{equation}
is not zero for some choice of data. Then we have a "Howe duality"
type of result. That is,\\
{\bf 3)} Suppose that integral \eqref{howe} is not zero for some
choice of data. Let $\nu$ be a finite nonarchimedean place where
all the data is unramified. Let $\theta_{\tau,\nu}, \epsilon_\nu,
\pi_\nu$ and $\omega_{\psi,\nu}$ denote the local representations
corresponding to the global ones which appear in integral
\eqref{howe}. Then, that integral induces an element in the space
$$Hom_{Sp_6\times
SO_{10}}((\theta_{\tau,\nu}\otimes\omega_{\psi,\nu})_U,\epsilon_\nu
\otimes\pi_\nu)$$ where $(\ldots)_U$ is the corresponding Jacquet
module. As in \cite{G2} section 6, it will be  proved in
\cite{G3}, that any two of the representations ${\tau_\nu},
\epsilon_\nu, \pi_\nu$ determines the third one uniquely.

From this we deduce that the representation $\pi_1$, as defined by
integral \eqref{lift1}, if nonzero, is the weak functorial lift
from $\epsilon_1$ and $\tau$. Similarly, if nonzero, integral
\eqref{lift2} defines an automorphic representation $\epsilon$ on
$Sp_6({\bf A})$, such that the representation $\pi$ is the weak
lift from $\epsilon$ and $\tau$. In the following Lemma, we will
prove that $\epsilon$, as defined by integral \eqref{lift2} is
cuspidal. We don't know if the image of the lift is irreducible,
however, it is clear that all irreducible summands are nearly
equivalent.

We start by proving\\
{\bf Lemma 6:} {\sl The representation $\epsilon$ is cuspidal.}\\
{\sl Proof:} Let $V$ be a standard unipotent subgroup of $Sp_6$.
By standard we mean that $V$ consists of upper unipotent matrices.
We shall use the form $\begin{pmatrix} &J_3\\ -J_3&
\end{pmatrix}$ to represent the symplectic group $Sp_6$ in terms of
matrices. Write
$$\widetilde{\theta}_\phi(l(u)(v,g))=\sum_{\xi\in Mat_{3\times
10}(F)}\omega_\psi(l(u)(v,g))\phi(\xi)=\sum_{\xi\in Mat_{3\times
10}(F)}\omega_\psi((X+\xi,Y,z)(v,g))\phi(0)$$ Here $\omega_\psi$
is the Weil representation and the above equalities are obtained
from the well known formulas for the Weil representation ( see for
example \cite{G-R-S3}).Plugging this identity into \eqref{lift2},
collapsing the summation and the integration over the $X$ variable
in $U$, we deduce that the integral $\int\limits_{V(F)\backslash
V({\bf A})}f(v)dv$ is zero for all choice of data, if and only if
the integral
\begin{equation}\label{lift3}
\int\limits_{SO_{10}(F)\backslash SO_{10}({\bf A})}\int
\limits_{V(F)\backslash V({\bf A})}\int \limits_{U_1(F)\backslash
U_1({\bf
A})}\varphi_\pi(g)\widetilde{\theta}_\tau(u_1(v,g))\psi_{1}(R_1)du_1dvdg
\end{equation}
is zero for all choice of data. Here $U_1$ is the subgroup of $U$
which consists of all matrices as in \eqref{mat} such that $X=0$.
The character $\psi_{1}$ is defined as follows. Write $u_1$ in the
coordinates as given in \eqref{mat}. If $R_1=(r_{i,j})\in
Mat_{3\times 3}$, then
$\psi_{1}(R_1)=\psi(r_{1,1}+r_{2,2}+r_{3,3})$.

At this point we need to consider the three nonconjugate maximal
unipotent radicals of $Sp_6$. We shall work out the details in the
case where $V$ is the unipotent radical of the Siegel parabolic.
The other two cases are treated similarly. Let
$$V=\left \{ v=\begin{pmatrix} I_3&Z\\ &I_3 \end{pmatrix} \ :\
Z\in Mat_{3\times 3}; Z^tJ_3=J_3Z \right \} $$ Thus \eqref{lift3}
is equal to
\begin{equation}\label{lift4}
\int\varphi_\pi(g)\widetilde{\theta}_\tau \left [ \begin{pmatrix}
I_3&Z&Y&R_1&R_2\\&I_3&&R_3&R_1^*\\
&&I_{10}&&Y^*\\ &&&I_3&Z^*\\ &&&&I_3 \end{pmatrix} (1,g)\right ]
\psi_{1}(R_1)du_1dZdg
\end{equation}
Here $g$ and $u_1$ are integrated as in \eqref{lift3} and $Z$ is
integrated over $V(F)\backslash V({\bf A})$. If $Z=(z_{i,j})$ we
consider the Fourier expansion of the above integral along the
variables $z_{2,3}, z_{3,2}$ and $z_{3,3}$ each integrated over
$F\backslash {\bf A}$. Conjugating by a suitable discrete matrix
in $R_3$ we collapse summation and integration, and deduce that
the vanishing of \eqref{lift4} for all choice of data is
equivalent to the vanishing of
\begin{equation}\label{lift5}
\int\varphi_\pi(g)\widetilde{\theta}_\tau \left [ \begin{pmatrix}
I_3&Z&Y&R_1&R_2\\&I_3&&&R_1^*\\
&&I_{10}&&Y^*\\ &&&I_3&Z^*\\ &&&&I_3 \end{pmatrix} (1,g)\right ]
\psi_{1}(R_1)du_1dZdg
\end{equation}
for all choice of data. Now $Z$ is integrated over $Mat_{3\times
3}(F)\backslash Mat_{3\times 3}({\bf A})$. Let $U_2$ denote the
unipotent subgroup of $SO_{22}$ which consists of all matrices of
the form
\begin{equation}\label{mat1}
\left \{ u=\begin{pmatrix} I_3&Z&Y&R_1&R_2\\&I_3&&&R_1^*\\
&&I_{10}&&Y^*\\ &&&I_3&Z^*\\ &&&&I_3 \end{pmatrix} \right \}
\end{equation}

Conjugating integral \eqref{lift5} by a suitable Weyl element, it
is then enough to show that the integral
\begin{equation}\label{lift6}
\int\varphi_\pi(g)\widetilde{\theta}_\tau \left [ \begin{pmatrix}
I_3&R_1&Y&Z&R_2\\&I_3&&&Z^*\\
&&I_{10}&&Y^*\\ &&&I_3&R_1^*\\ &&&&I_3 \end{pmatrix} (1,g)\right ]
\psi_{1}(R_1)du_2dg
\end{equation}
is zero for all choice of data. Here $u_2$ is integrated over
$U_2(F)\backslash U_2({\bf A})$. Next we consider the Fourier
expansion of \eqref{lift6} along the unipotent group
$$\begin{pmatrix} I_3&&&&\\&I_3&&R_3&\\
&&I_{10}&&\\ &&&I_3&\\ &&&&I_3 \end{pmatrix}\ \ \ \ R_3\in
Mat^0_{3\times 3}$$ with points in $F\backslash {\bf A}$. Under
the action of $GL_3(F)$ embedded as $diag(m,m,I_{10},m^*,m^*)$,
there are two orbits to consider in this expansion. The first one
is given by a sum of integrals of the form
\begin{equation}\label{lift7}
\int\varphi_\pi(g)\widetilde{\theta}_\tau \left [ \begin{pmatrix}
I_3&R_1&Y&Z&R_2\\&I_3&&R_3&Z^*\\
&&I_{10}&&Y^*\\ &&&I_3&R_1^*\\ &&&&I_3 \end{pmatrix} (1,g)\right ]
\psi_{1}(R_1)\psi_2(R_3)du_2dR_3dg
\end{equation}
where $\psi_2(R_3)=\psi(r_{1,2})$ for all $R_3=(r_{i,j})$. After
some further tedious expansions and suitable conjugations one can
show that the above integral is a sum of Fourier coefficients
where each one of them corresponds to a unipotent orbit which is
bigger than or equal to $(41^{18})$. Since ${\cal
O}_{SO_{22}}(\widetilde{\theta}_\tau)=(3^71)$ it follows that
\eqref{lift7} is zero for all choice of data. Thus we are left
with the trivial orbit. That is, we need to show that for all
choice of data the integral
\begin{equation}\label{lift8}
\int\varphi_\pi(g)\widetilde{\theta}_\tau \left [ \begin{pmatrix}
I_3&R_1&Y&Z&R_2\\&I_3&&R_3&Z^*\\
&&I_{10}&&Y^*\\ &&&I_3&R_1^*\\ &&&&I_3 \end{pmatrix} (1,g)\right ]
\psi_{1}(R_1)du_2dR_3dg
\end{equation}
is zero. Next we expand integral \eqref{lift8} along the unipotent
set
$$\begin{pmatrix} I_3&&&&\\&I_3&X&&\\
&&I_{10}&X^*&\\ &&&I_3&\\ &&&&I_3 \end{pmatrix}\ \ \ \ X\in
Mat_{3\times 10}$$ with points in $F\backslash {\bf A}$. The group
$GL_3(F)\times SO_{10}(F)$ acts on this expansion with various
orbits. The orbits can be identified with vectors $\xi_i\in
F^{10}$, for $1\le i\le 3$, according to the length of these
vectors. If at least one of these vectors has nonzero length, then
the corresponding Fourier coefficient contains as inner
integration a Fourier coefficient which  corresponds to the
unipotent orbit $(51^{17})$. By the smallness of
$\widetilde{\theta}_\tau$ this Fourier coefficient is zero. Thus
we are left with the cases where the length of all vectors is
zero. Hence we can write the above Fourier expansion as a sum of
integrals of the form
\begin{equation}\label{lift9}
\int\varphi_\pi(g)\widetilde{\theta}_\tau \left [ \begin{pmatrix}
I_3&R_1&Y&Z&R_2\\&I_3&X&R_3&Z^*\\
&&I_{10}&X^*&Y^*\\ &&&I_3&R_1^*\\ &&&&I_3 \end{pmatrix}
(1,g)\right ] \psi_{1}(R_1)\psi_\nu(X)dXdu_2dR_3dg
\end{equation}
Here $\nu=(\nu_1,\nu_2,\nu_3)$ is one of the vectors $(0,0,0),
(1,0,0), (1,1,0)$ or $(1,1,1)$. The characters $\psi_\nu$ are
defined as follows. Given $X=(x_{i,j})\in Mat_{3\times 10}$ we
define $\psi_\nu(X)=\psi(\nu_1 x_{1,1}+\nu_2 x_{2,2}+\nu_3
x_{3,3})$ where $\nu$ is any one of the above four vectors. The
variable $g$ is integrated over
$SO_{10-2(\nu_1+\nu_2+\nu_3)}(F)U_\nu\backslash SO_{10}({\bf A})$,
where $U_\nu$ is the unipotent radical of the maximal parabolic
subgroup of $SO_{10}$ whose Levi part is
$GL_{\nu_1+\nu_2+\nu_3}\times SO_{10-2(\nu_1+\nu_2+\nu_3)}$. If
$\nu=(0,0,0)$ then from the definition of
$\widetilde{\theta}_\tau$ we obtain the integral
$\int\varphi_\pi(g)\theta_\tau(g)dg$ as inner integration. Here
$g$ is integrated over $SO_{10}(F)\backslash SO_{10}({\bf A})$.
Hence by the cuspidality of $\pi$ this integral is zero. If $\nu$
is any one of the other three cases, we continue by expanding
integral \eqref{lift9} along the group $U_\nu(F)\backslash
U_\nu({\bf A})$. If this group is not abelian, we first expand
along its center. By the smallness of $\widetilde{\theta}_\tau$ we
obtain that all orbits in the various expansions are zero except
the constant term along $U_\nu$. This follows from the fact that
each nontrivial Fourier coefficient will contain, as inner
integration, a Fourier coefficient corresponding to a unipotent
orbit which is greater than or not related to $(3^71)$. Hence it
will give zero contribution. The constant term will also vanish.
Indeed, in this case after factoring the integral, we will obtain
the integral $\int \varphi_\pi(ug)du$ as inner integration. Here,
the variable $u$ is integrated over $U_\nu(F)\backslash U_\nu({\bf
A})$. Hence we get zero by the cuspidality of $\pi$. This
completes the proof of the vanishing for the unipotent radical of
the Siegel parabolic. We still have to prove the same for the
other two maximal unipotent radicals. This is done in a similar
way and we shall omit the details.   \hfill $\blacksquare$

\section {\bf The Main Theorem}

We recall the basic notations we used in the previous sections.
Let $\pi$ denote a generic cuspidal representation of
$GSO_{10}({\bf A})$. We shall assume that $\pi$ has a trivial
central character. Let $\tau$ denote a cuspidal representation of
$PGL_2({\bf A})$. In section two we introduced the global integral
\begin{equation}\label{main1}
\int\limits_{Z({\bf A})G(F)\backslash G({\bf A})}\varphi_\pi(g)
\theta_\tau(g)E(g,s)dg
\end{equation}
Here $G=GSO_{10}$, and $Z$ is the center of $G$. The
representation $\theta_\tau$ was introduced in section three, and
the Eisenstein series $E(g,s)$ was defined at the beginning of
section two. It follows from sections 2 and 3 that \eqref{main1}
unfolds to a Whittaker integral, and hence is Eulerian. From
Proposition 4 in section 4 we deduce that at the
non-archimedean unramified places this integral represents the $L$
function $L(Spin_{10}\times St,\pi\times \tau, 2s-1/2)$. Let $S$
denote a finite number of places in $F$ which includes all
archimedean places, such that outside of $S$ all data is
unramified. We shall denote by $L^S(Spin_{10}\times St,\pi\times
\tau, 2s-1/2)$ the corresponding partial $L$ function.

It follows from \cite{K-R} that the Eisenstein series $E(g,s)$ can
have at most a simple pole at the points $s=1$ and $s=3/4$. The
residue at $s=1$ of $E(g,s)$ is the trivial representation, and
hence the residue of integral \eqref{main1} is zero at that point.
We are interested in the residue of \eqref{main1} at the other
point. In this section only, we shall denote by $\theta$ the
residual representation of   $E(g,s)$ at $s=3/4$. In other words,
we denote $\theta(g)=Res_{s=3/4}E(g,s)$. It follows from
\cite{G-R-S4} section 6 that $\theta$ is the  irreducible minimal
representation of $G$. We are now ready to state and prove our
main result\\
{\bf Main Theorem:} {\sl Let $\pi$ be an irreducible generic
cuspidal representation of the group $G(\bf A)$ which has a
trivial central character. Then the following are equivalent:\\
1) The partial $L$ function $L^S(Spin_{10}\times St,\pi\times
\tau, 2s-1/2)$ has a simple pole at $s=3/4$.\\
2) The period integral
\begin{equation}\label{main2}
\int\limits_{SO_{10}(F)\backslash SO_{10}({\bf A})}\varphi_\pi(g)
\theta_\tau(g)\theta(g)dg
\end{equation}
is nonzero for some choice of data.\\
3) There exists a generic cuspidal representation $\sigma$ of the
exceptional group $G_2({\bf A})$ such that $\pi$ is the weak lift
from the representation $\sigma\times\tau$ of the group $G_2({\bf
A})\times PGL_2({\bf A})$.
}\\
{\sl Proof:} Two parts are not hard to prove. Suppose 1) holds.
Arguing as in \cite{Ga-S} section 7, one can show that for any
place in the global field $F$, given a value of $s\in {\bf C}$,
one can choose data such that the local integral \eqref{local1} is
nonzero. From the results in the previous sections we hence deduce
that the integral
\begin{equation}\label{main3}
\int\limits_{Z({\bf A})GSO_{10}(F)\backslash GSO_{10}({\bf
A})}\varphi_\pi(g) \theta_\tau(g)\theta(g)dg
\end{equation}
is nonzero for some choice of data. Factoring the similitude in
\eqref{main3}, we deduce that if the partial $L$ has a pole at
$s=3/4$, then integral \eqref{main2} is nonzero for some choice of
data. Thus 1) implies 2).

The implication 3) implies 1) is also not hard. Suppose that $\pi$
is the weak  lift from $\sigma\times \tau$. Then branching from
$SO_{10}({\bf C})$ to $G_2({\bf C})\times SL_2({\bf C})$ we obtain
that $L^S(Spin_{10}\times St,\pi\times\tau, 2s-1/2)$ is equal to
$$L^S(G_2\times Sym^2,\sigma\times\tau,2s-1/2)
L^S(G_2,\sigma,2s-1/2)L^S(Sym^2,\tau,2s-1/2)\zeta^S(2s-1/2)$$ Here
$G_2$ is the standard seven dimensional representation of
$G_2({\bf C})$ and $Sym^2$ is the symmetric square representation
of $SL_2({\bf C})$. Also, $\zeta^S(2s-1/2)$ denotes the partial
global $\zeta$ function. Since $\sigma$ is generic it follows from
\cite{G-R-S5} that $\sigma$ has a nontrivial lift to a generic
cuspidal representation of $Sp_6({\bf A})$ or to a cuspidal
representation on $PGL_3({\bf A})$. Using the lifting of generic
cuspidal representations from Classical groups to $GL_n$ as proved
in \cite{C-K-PS-S},  we deduce that the above product of $L$
functions is actually a product of partial $L$ functions
associated with tensor product representations on certain $GL$'s.
Hence, at $s=3/4$ these $L$ functions cannot vanish. On the other
hand, $\zeta^S(2s-1/2)$ has a simple pole at that point. Hence 3)
implies 1).

Finally we prove that 2) implies 3). Before providing the details,
let us first explain the idea of the proof. Assuming \eqref{main2}
is nonzero for some choice of data, we need to construct a generic
cuspidal representation $\sigma$ of the group $G_2({\bf A})$ and
prove that $\pi$ is the weak lift of $\sigma\times\tau$. Suppose
that there exists a cuspidal representation $\epsilon$ of
$Sp_6({\bf A})$ such that the integral
\begin{equation}\label{sp1}
f^{R,\psi}(h)=\int\limits_{SL_2(F)\backslash SL_2({\bf A})}
\int\limits_{V(F)\backslash V({\bf A})}f(vm)\psi_V(v)dvdm
\end{equation}
is nonzero for some choice of data. Here $f$ is a vector in the
space of $\epsilon$ and the group $R=V\cdot SL_2$ is defined as
follows. Let $V$ denote the unipotent radical of the maximal
parabolic subgroup of $Sp_6$ whose levi part is $GL_2\times SL_2$.
Thus in matrices we have
$$ V\cdot SL_2=\left \{ \begin{pmatrix} I_2&x&y\\ &I_2&x^*\\ &&I_2
\end{pmatrix}\begin{pmatrix} m&&\\ &m& \\ &&m \end{pmatrix} :
x\in Mat_{2\times 2}; y\in Mat'_{2\times 2}; m\in SL_2 \right \}
$$ Here $Mat'_{2\times 2}=\{ y\in Mat_{2\times 2} : y^tJ_2=J_2y
\}$. Also, we define $\psi_V(v)=\psi(trx)$. It is also not hard to
prove that if $\epsilon$ is such that integral \eqref{sp1} is
nonzero for some choice of data, then $\epsilon$ must be generic.
Assuming all this, it follows from \cite{G-Ja} ( see also
\cite{G-H} the discussion before and after equation (27)) that
$\epsilon$ is the weak functorial lift from a cuspidal
representation $\sigma$ on $G_2({\bf A})$.

Let $\epsilon$ denote the representation defined by right
translations of all functions given by integral \eqref{lift2}. As
proved in section 5 we know that $\epsilon$ is a cuspidal
representation of $Sp_6({\bf A})$. Thus if we can prove that
$f^{R,\psi}(h)$ is nonzero for some choice of function $f(h)$ as
in \eqref{lift2} it will follow from section 5 that $\pi$ is a
weak endoscopic lift from $\epsilon\times \tau$. From the above
discussion it will follow that $\epsilon$ is a lift from a generic
cuspidal representation $\sigma$ of $G_2({\bf A})$. This will
prove that 2) implies 3).

Thus, to conclude the proof of the theorem we will show that if
\eqref{main2} is not zero for some choice of data, then the
integral
\begin{equation}\label{main4}
\int\limits_ {SL_2(F)\backslash SL_2({\bf A})}
\int\limits_{V(F)\backslash V({\bf A})}\int\int
\limits_{U(F)\backslash U({\bf A})}\varphi_\pi(g)
\widetilde{\theta}_\phi(l(u)(vm,g))\widetilde{\theta}_\tau(u(vm,g))
\psi_V(v)dudgdvdm
\end{equation}
is not zero for some choice of data. Here $g$ is integrated over
$SO_{10}(F)\backslash SO_{10}({\bf A})$. Since the inner
integrations, over the variables $u$ and $g$ produce a cusp form
on $Sp_6({\bf A})$, the integral over $SL_2(F)\backslash SL_2({\bf
A})$ converges. At some point we will need to interchange the
order of the integration. The justification for this is given in
\cite{K-R} Proposition 5.3.1 . In other words, one needs to add
the action of the Casimir element for $SL_2$ inside the Schwartz
function inside the theta representation. Doing that, the theta
representation becomes rapidly decreasing as a function of the
$SL_2$. Because of \cite{G-R-S4} Theorem 6.9 this will not harm
the generality of our argument.

We start by unfolding the theta series. As in section 5 we have
$$\widetilde{\theta}_\phi(l(u)(vm,g))=\sum_{\eta,\xi}\omega_\psi(
l(u)(vm,g))\phi(\eta,\xi)=\sum_{\eta,\xi}\omega_\psi(
l(\eta)l(u)(vm,g))\phi(0,\xi).$$ Here $\eta\in Mat_{2\times 10}$
and $\xi\in F^{10}$. Plugging this into \eqref{main4} we collapse
the summation over $\eta$ with the suitable integration inside
$U$. We then consider the Fourier expansion along the unipotent
subgroup of $SO_{22}$ defined by $I_{22}+r_1e'_{1,3}+r_2e'_{1,4}+
r_2e'_{2,3}+r_2e'_{2,4}+r_5e'_{2,6}$. Here
$e'_{i,j}=e_{i,j}-e_{23-j,23-i}$ and $e_{i,j}$ is the matrix which
has one at the $(i,j)$ position and zero elsewhere. Once again,
collapsing summation and integration, integral \eqref{main4}
equals
\begin{equation}\label{main5}
\int
\varphi_\pi(g)\sum_{\xi}\omega_\psi(l(u_2)(m,g)l(u_1))\phi(0,\xi)
\widetilde{\theta}_\tau(u_3u_2(m,g)u_1)\psi_{U_3}(u_3)d(...)
\end{equation}
Here, the variables $m$ and $g$ are integrated as in
\eqref{main4}. In term of matrices we have
\begin{equation}\label{ma1}
u_1=\begin{pmatrix} I_2&&&&&&\\&I_2&&&z_1&&\\&&I_2&z_3&z_2&z_1^*&\\
&&&I_{10}&z_3^*&&\\ &&&&I_2&&\\&&&&&I_2&\\&&&&&&I_2 \end{pmatrix}
u_2=\begin{pmatrix} I_2&&&&&&\\&I_2&&y_1&&y_2&&\\&&I_2&&&&\\
&&&I_{10}&&y_1^*&\\ &&&&I_2&&\\&&&&&I_2&\\&&&&&&I_2
\end{pmatrix}
\end{equation}
Here, the variable $z_1$ is integrated over $Mat_{2\times 2}({\bf
A}) $, the variable $z_2$ is integrated over $Mat_{2\times
2}^0({\bf A})$ ( see the definition right after \eqref{mat}) and
$z_3$ is integrated over $Mat_{2\times 10}({\bf A}) $. The
variable $y_1$ is integrated over $Mat_{2\times 10}(F)\backslash
Mat_{2\times 10}({\bf A}) $ and $y_2$ is integrated over
$Mat_{2\times 2 }^0(F)\backslash Mat_{2\times 2}^0({\bf A})$. The
variable $u_3$ is defined as
\begin{equation}\label{mat2}
u_3=u_3(x_1,x_2,x_3,r_1,r_2,r_3,r_4)=\begin{pmatrix} I_2&x_1&x_2&r_1&r_2&r_3&r_4\\
&I_2&x_3&&&&r_3^*\\&&I_2&&&&r_2^*\\
&&&I_{10}&&&r_1^*\\ &&&&I_2&x_3^*&x_2^*\\&&&&&I_2&x_1^*\\&&&&&&I_2
\end{pmatrix}
\end{equation}
Here all variables are integrated over $Mat_{a\times b}$, for
suitable choice of $a$ and $b$  with points in $F\backslash {\bf
A}$, except $r_4$ which is integrated over $Mat^0_{2\times 2}$.
Finally, the character $\psi_{U_3}$ is defined as follows. For
$r_2=(r_{i,j})$ and $x_3=(x_{i,j})$ as described in \eqref{mat2}
we set $\psi_{U_3}(u_3)=\psi(r_{1,1}+r_{2,2}+x_{1,1}+x_{2,2})$.
Notice that the variable $v$ which appeared inside the Weil
representation one equation before \eqref{main5} does not appear
inside $\omega_\psi$ in \eqref{main5}. This follows from the fact
that $\omega_\psi((v,1))\phi(0,\xi)=\phi(0,\xi)$.

The group $u_3(0,r_4)$ is the one dimensional unipotent group
which is the center of the maximal unipotent subgroup of
$SO_{22}$. Let $N_2$ denote the standard unipotent radical of the
maximal parabolic subgroup of $SO_{22}$ whose levi part is
$GL_2\times SO_{18}$. We have
$$\int\limits_{F\backslash {\bf A}}\widetilde
{\theta}_\tau(u_3(0,r_4))dr_4=\sum_{\alpha}
\int\limits_{N_2(F)\backslash N_2({\bf A})}
\widetilde{\theta}_\tau(n_2)\psi(\alpha\cdot n_2)dn_2$$ where the
sum is over all $\alpha\in Mat_{2\times 18}(F)$. The group
$SO_{18}(F)$ acts on this expansion. If $\alpha=\begin{pmatrix}
\alpha_1\\ \alpha_2 \end{pmatrix}$ where $\alpha_i\in F^{18}$,
then the various orbits can be parameterized by the rank of the
matrix $\alpha$ and by the length of $\alpha_i$.

Plugging this expansion into \eqref{main5} we claim that only one
orbit contributes a nonzero term. Indeed, after plugging in
\eqref{main5}, we conjugate the matrix
$u_3(x_1,x_2,0,r_1,r_2,r_3,0)$ from right to left. Because of the
character $\psi_{U_3}$ we get zero contribution unless the rank of
$\alpha$ is two and both vectors $\alpha_i$ are nonzero and have
zero length. From this we deduce that \eqref{main5} equals
\begin{equation}\label{main6}
\int
\varphi_\pi(g)\sum_{\xi}\omega_\psi(l(u_2)(m,g)l(u_1))\phi(0,\xi)
\sum_{\gamma}\widetilde{\theta}_\tau^{N_2,\psi} (\gamma
u_3u_2(m,g)u_1)\psi_{U_3}(u_3)d(...)
\end{equation}
Here $\gamma\in P^0_2(SO_{14})(F)\backslash SO_{18}(F)$ where
$P^0_2(SO_{14})$ is the subgroup of $SO_{18}$ which consists of
all matrices of the form
$$\begin{pmatrix} I_2&*&*\\&k&*\\&&I_2 \end{pmatrix} \ \ \ k\in SO_{14}$$
Also, we denote
$$\widetilde{\theta}_\tau^{N_2,\psi}(k)=
\int\limits_{N_2(F)\backslash N_2({\bf A})}
\widetilde{\theta}_\tau(n_2k)\psi_{N_2}( n_2)dn_2\ \ \ \ \ \ \
k\in SO_{22}$$ where $\psi_{N_2}$ is defined as follows.
Identifying $N_2$ with the group of all matrices of the form
$n_2=u_3(x_1,x_2,0,r_1,r_2,r_3,r_4)$, we set
$\psi_{N_2}(n_2)=\psi(trx_1)$.

Next we consider the double coset space $P^0_2(SO_{14})\backslash
SO_{18}/P_2(SO_{14})$. Here $P_2(SO_{14})$ is the standard
parabolic subgroup of $SO_{18}$ whose levi part is $GL_2\times
SO_{14}$. Thus $P_2(SO_{14})$ contains $P^0_2(SO_{14})$.
Considering the various orbits and arguing as above, one can show
that only one orbit, which we shall describe bellow, contributes a
nonzero term to \eqref{main6}. Let
$$w_1=\begin{pmatrix} I_2&&&&&&\\&&I_2&&&& \\&I_2&&&&&
\\&&&I_{10}&&& \\&&&&&I_2& \\&&&&I_2&& \\&&&&&&I_2 \end{pmatrix}\
\ \ \ \ t(\delta)=\begin{pmatrix} I_2&&&&&&\\&I_2&\delta &&&&&\\&&I_2&&&&\\
&&&I_{10}&&&\\ &&&&I_2&\delta^* &\\&&&&&I_2&\\&&&&&&I_2
\end{pmatrix}
$$ where $\delta\in Mat_{2\times 2}$. Then integral \eqref{main6}
equals
\begin{equation}\label{main7}
\int
\varphi_\pi(g)\sum_{\xi}\omega_\psi(l(u_2)(m,g)l(u_1))\phi(0,\xi)
\sum_{\gamma,\delta}\widetilde{\theta}_\tau^{N_2,\psi}
(w_1t(\delta)\gamma u_3u_2(m,g)u_1)\psi_{U_3}(u_3)d(...)
\end{equation}
Here $\gamma\in P^0_2(SO_{10})(F)\backslash SO_{14}(F)$ where
$P^0_2(SO_{10})$ is the subgroup of $SO_{14}$ which consists of
all matrices in of the form
$$\begin{pmatrix} I_2&*&*\\&k&*\\&&I_2 \end{pmatrix} \ \ \ k\in SO_{10}$$
Also, we have $\delta\in Mat_{2\times 2}(F).$ In \eqref{main7} we
conjugate the matrix $u_3(x_1,0)$ across $t(\delta)$. We obtain
from the commutation relations the matrix $u_3(0,x_1\delta,0)$.
When we conjugate this matrix across $w_1$ and change variables in
$N_2$, we obtain the integral
$$\int\limits_{(F\backslash {\bf A})^4}\psi(tr(x_1\delta))dx_1$$
as inner integration. Clearly this integral is zero unless
$\delta=0$. Thus, if $\delta\ne 0$ in \eqref{main7} we get zero
contribution.

Next we consider the space $P^0_2(SO_{10})\backslash
SO_{14}/P_2(SO_{10})$. Here $P_2(SO_{10})$ is the standard
parabolic subgroup of $SO_{14}$ whose levi part is $GL_2\times
SO_{10}$. Thus $P_2(SO_{10})$ contains $P^0_2(SO_{10})$. As before
one can check that all orbits contribute zero except one. Let
$$w_2=\begin{pmatrix} I_2&&&&&&\\&I_2&&&&& \\&&&&I_2&&
\\&&&I_{10}&&& \\&&I_2&&&& \\&&&&&I_2& \\&&&&&&I_2 \end{pmatrix}\
\ \ \ \ t(\delta,\mu)=\begin{pmatrix} I_2&&&&&&\\&I_2&&&&&&\\&&I_2&\delta &\mu&&\\
&&&I_{10}&\delta^*&&\\ &&&&I_2&&\\&&&&&I_2&\\&&&&&&I_2
\end{pmatrix}
$$ where $\delta\in Mat_{2\times 10}$ and $\mu\in Mat_{2\times
2}^0$. Conjugating by the matrices $u_3(0,x_2,0,r_1,0)$ across
from left to right, as we did before with $u_3(x_1,0)$, we obtain
that only when $\delta=0$ and $\mu=0$ then we get a nonzero
contribution. Thus, integral \eqref{main7} equals
\begin{equation}\label{main8}
\int
\varphi_\pi(g)\sum_{\xi}\omega_\psi(l(u_2)(m,g)l(u_1))\phi(0,\xi)
\widetilde{\theta}_\tau^{N_2,\psi} (w_1w_2
u_3(x_3)u_2(m,g)u_1)\psi(trx_3)d(...)
\end{equation}
where we wrote $u_3(x_3)$ for $u_3(0,0,x_3,0)$.

Next, we consider the Fourier expansion
$$\widetilde{\theta}_\tau^{N_2,\psi}(k)=\sum_{\beta\in
F}\int\limits_{F\backslash {\bf
A}}\widetilde{\theta}_\tau^{N_2,\psi}((I_{22}+
re'_{3,19})k)\psi(\beta r)dr$$ One can check that if $\beta\ne 0$
then the corresponding summand is zero. Indeed, in this case we
obtain as inner integration, Fourier coefficients of
$\widetilde{\theta}_\tau$ corresponding to  unipotent orbits which
are greater or equal to $(41^{18})$. By section 5 point ${\bf b)}$
it follows that all these Fourier coefficients are zero. Thus we
are left with the summand corresponding to $\beta=0$. This we can
further expand along the set
$$\begin{pmatrix} I_2&&&\\&I_2&x&*&\\
&&I_{12}&x^*&&\\ &&&I_2&\\&&&&I_2
\end{pmatrix}$$
where $x\in Mat_{2\times 12}.$ The group $SO_{12}$ acts on this
expansion, and as before we can parameterize the various orbits by
matrices $\alpha=\begin{pmatrix} \alpha_1\\
\alpha_2\end{pmatrix}\in Mat_{2\times 12}(F)$ where the invariants
are the rank of $\alpha$ and the length of $\alpha_i$. Once again
we are left with the contribution of one orbit which corresponds
to rank two matrices such that both $\alpha_i$ are nonzero vectors
of length zero. Thus, integral \eqref{main8} is equal to
\begin{equation}\label{main9}
\int
\varphi_\pi(g)\sum_{\xi}\omega_\psi(l(u_2)(m,g)l(u_1))\phi(0,\xi)
\sum_{\gamma}\widetilde{\theta}_\tau^{N_4,\psi} (\gamma w_1w_2
u_3(x_3)u_2(m,g)u_1)\psi(trx_3)d(...)
\end{equation}
Here, the group $N_4$ consists of all unipotent matrices in
$SO_{22}$ of the form
$$\begin{pmatrix} I_2&r_1&*&*&*&*&*\\&I_2&r_2&*&*&*&*&\\&&I_2&&&*&*\\
&&&I_{10}&&*&*\\ &&&&I_2&r_2^*&*\\&&&&&I_2&r_1^*\\&&&&&&I_2
\end{pmatrix}$$ and
$$\widetilde{\theta}_\tau^{N_4,\psi}(k)=
\int\limits_{N_4(F)\backslash N_4({\bf A})}
\widetilde{\theta}_\tau(n_4k)\psi_{N_4}( n_4)dn_4$$ where
$\psi_{N_4}(n_4)=\psi(tr(r_1+r_2)).$ Also, we have $\gamma\in
P^0_2(SO_{10})(F)\backslash SO_{14}(F).$ Continuing this process
of Fourier expansions it follows from the smallness properties of
$\widetilde{\theta}_\tau$ that
$\widetilde{\theta}_\tau^{N_4,\psi}(k)=\widetilde{\theta}_\tau^{N_6,\psi}(k)$
where
$$\widetilde{\theta}_\tau^{N_6,\psi}(k)=
\int\limits_{N_6(F)\backslash N_6({\bf A})}
\widetilde{\theta}_\tau(n_6k)\psi_{N_6}( n_6)dn_6$$ Here $N_6$ is
the  unipotent subgroup of $SO_{22}$ which consists of all
matrices of the form
$$\begin{pmatrix} I_2&r_1&*&*&*&*&*\\&I_2&r_2&*&*&*&*&\\&&I_2&r_3&*&*&*\\
&&&I_{10}&r_3^*&*&*\\ &&&&I_2&r_2^*&*\\&&&&&I_2&r_1^*\\&&&&&&I_2
\end{pmatrix}$$ and $\psi_{N_6}$ is the character $\psi_{N_4}$
extended trivially from $N_4$ to $N_6$.

As before we consider the space $P^0_2(SO_{10})\backslash
SO_{14}/P_2(SO_{10})$. We obtain a nonzero contribution from one
term. Thus, integral \eqref{main9} equals
\begin{equation}\label{main10}
\int
\varphi_\pi(g)\sum_{\xi}\omega_\psi(l(u_2)(m,g)l(u_1))\phi(0,\xi)
\sum_{\delta,\mu}\widetilde{\theta}_\tau^{N_6,\psi}
(w_2t(\delta,\mu) w_1w_2 u_2(m,g)u_1)d(...)
\end{equation}
Conjugating $w_2w_1t(\delta,\mu) w_1w_2$ we see that as a group,
this group of matrices coincides with the group $U_2$. Thus we may
collapse summation with integration. Hence \eqref{main10} is equal
to
\begin{equation}\label{main11}
\int\int\limits_{ (U_1U_2)({\bf A})}
\varphi_\pi(g)\sum_{\xi}\omega_\psi((m,g)l(u_2u_1))\phi(0,\xi)
\widetilde{\theta}_\tau^{N_6,\psi}
(w_2w_1w_2(m,g)u_2u_1)du_1du_2dgdm
\end{equation}
Here the variables $m$ and $g$ are integrated as before. That is,
the variable $m$ is integrated over $SL_2(F)\backslash SL_2({\bf
A})$, and $g$ is integrated over $SO_{10}(F)\backslash
SO_{10}({\bf A})$.

Notice that $N_6=UV_6$ where $U$ is the unipotent radical of the
parabolic subgroup $Q$ as defined in \eqref{mat}. The group $V_6$
is the group of all unipotent matrices of the form
$$V_6=\left \{ \begin{pmatrix} I_2&r_1&r_2&&&&\\&I_2&r_3&&&& \\
&&I_2&&&&\\ &&&I_{10}&&& \\ &&&&I_2&r_3^*&r_2^*\\ &&&&&I_2&r_1^*\\
&&&&&&I_2
\end{pmatrix} \right \} $$

To complete the proof of the theorem, we assume that part 2) in
the statement of the Theorem holds, and assume that integral
\eqref{main11} is zero for all choice of data. We shall derive a
contradiction. First, using arguments similar as to those in
\cite{Ga-S} section 7,  we may deduce that if indeed
\eqref{main11} is zero for all choice of data, then so is the integral
\begin{equation}\label{main12}
\int\limits_{SL_2(F)\backslash SL_2({\bf A})}\int\limits_{
SO_{10}(F)\backslash SO_{10}({\bf A})}
\varphi_\pi(g)\sum_{\xi}\omega_\psi((m,g))\phi_2(\xi)
\widetilde{\theta}_\tau^{N_6,\psi} (w_2w_1w_2(m,g))dgdm
\end{equation}
Here we wrote $\phi=\phi_1\otimes\phi_2$ where $\phi_1$ is a
Schwartz function of $Mat_{2\times 10}({\bf A})$ and $\phi_2$ is a
Schwartz function of ${\bf A}^{10}$. Since $w_2w_1w_2$ normalizes
the group of matrices $(m,g)$ as above,  we conjugate it to the
right and may ignore it.

Next, we claim that $\widetilde{\theta}_\tau^{N_6,\psi}((m,g))=
\widetilde{\theta}_\tau^{N_6,\psi}((1,g))$ for all $m\in SL_2({\bf
A})$ embedded in $SO_{22}({\bf A})$ as in \eqref{main12}. We will
prove it in Lemma 7 below. Assuming that, we deduce that if
\eqref{main12} is zero for all choice of data, then so is the
integral
\begin{equation}\label{main13}
\int\limits_{ SO_{10}(F)\backslash SO_{10}({\bf A})}
\varphi_\pi(g)f(g)\widetilde{\theta}_\tau^{N_6,\psi} ((1,g))dg
\end{equation}
where we denoted
$$f(g)=\int\limits_{SL_2(F)\backslash SL_2({\bf
A})}\sum_{\xi}\omega_\psi((m,g))\phi_2(\xi)dm$$ Notice that the
inner summation is in fact the theta representation defined on the
group $\widetilde{Sp}_{20}({\bf A})$. As explained in the
beginning of the proof, after possible adjustments of the Schwartz
function, it follows from \cite{G-R-S4} Theorem 6.9, that $f(g)$
is actually the minimal representation of $SO_{10}$ which we
denoted by $\theta$. Hence, from \eqref{main13} and the above
discussion, we may assume that the integral
\begin{equation}\label{main14}
\int\limits_{ SO_{10}(F)\backslash SO_{10}({\bf A})}
\varphi_\pi(g)\theta(g)\widetilde{\theta}_\tau^{N_6,\psi}
((1,g))dg
\end{equation}
is zero for all choice of data. Recall that $N_6=UV_6$. Also the
character $\psi_{N_6}$ is trivial on $U_6$. In Lemma 7 below, we
shall prove that
\begin{equation}\label{speh1}
\int\limits_{V_6(F)\backslash V_6({\bf A})}
E_\tau(vh)\psi_{V_6}(v)dv
\end{equation}
is nonzero for some choice of data. Here $E_\tau$ is the Speh
representation as was defined at the beginning of section 5, and
by abuse of notations, we view $V_6$ as a subgroup of $GL_6$. The
character $\psi_{V_6}$ is the restriction of $\psi_{N_4}$ to the
group $V_6$.

As a function of $g\in SO_{10}({\bf A})$, the constant term
$\widetilde{\theta}_\tau^{U}((1,g))$, is realized in the
representation space of $\theta_\tau$. In other words, in view of
the nonvanishing property of \eqref{speh1}, it follows that the
vanishing for all data of \eqref{main14} implies that
\eqref{main2} is zero for all choice of data. But this contradicts
our assumption. Thus 2) implies 3). This completes the proof of
the main theorem.   \hfill $\blacksquare$

We still need to prove\\
{\bf Lemma 7:} {\sl Let $E_\tau$ denote the Speh representation on
$GL_6({\bf A})$ as defined at the beginning of section 5. Then the
integral
\begin{equation}\label{speh2}
{\cal F}(h)=\int\limits_{(Mat_{2\times 2}(F)\backslash
Mat_{2\times 2}({\bf A}))^3}E_\tau\left ( \begin{pmatrix}
I_2&x&z\\&I_2&y\\&&I_2
\end{pmatrix} h\right ) \psi(tr(x+y))dxdydz
\end{equation}
is nonzero for some choice of data. Moreover, if $h=diag(m,m,m)$,
where $m\in SL_2({\bf A})$, then the ${\cal F}(h)={\cal F}(e)$. }\\
{\sl Proof:} It follows from \cite{G4} Proposition 5.3 that ${\cal
O}(E_\tau)=(3^2)$. This means that $E_\tau$ has no nonzero Fourier
coefficients for any unipotent orbit which is bigger than or not
related to $(3^2)$. Let $t(r)=diag(\begin{pmatrix} 1&r\\ &1
\end{pmatrix},\begin{pmatrix} 1&r\\ &1 \end{pmatrix},
\begin{pmatrix} 1&r\\ &1 \end{pmatrix})$. We expand ${\cal F}(h)$
along the group of matrices $t(r)$ with points in $F\backslash
{\bf A}$. We have
$${\cal F}(h)=\sum_{\alpha\in F}\int\limits_{F\backslash {\bf A}}
{\cal F}(t(r)h)\psi(\alpha r)dr$$ If $\alpha$ is not zero then the
corresponding expansion consists of  Fourier coefficients for
$E_\tau$ which correspond to unipotent orbits which are greater
than or not related to $(3^2)$. Hence, if $\alpha\ne 0$, we get
zero contribution which means that ${\cal F}(t(r)h)={\cal F}(h)$
for all $r\in {\bf A}$. Let $\nu=\begin{pmatrix} 0&1\\ 1&0
\end{pmatrix}$. Denote $w=diag(\nu,\nu,\nu)$. Then clearly
${\cal F}(wh)={\cal F}(h)$. Since $w$ and $t(r)$ generate
$SL_2({\bf A})$, the second assertion follows.

As for the first part, we shall assume that ${\cal F}(h)$ is zero
for all choice of data and derive a contradiction. Let $w$ denote
the Weyl element of $GL_6$ with one at positions $(1,1); (2,3);
(3,5); (4,2); (5,4); (6,6)$ and zero elsewhere. Conjugating $w$
from left to right we deduce that the integral
\begin{equation}\label{speh3}
\int\limits_{(F\backslash {\bf A})^{12}}E_\tau\left (
\begin{pmatrix} 1&x_1&x_2&&r_1&r_2\\ &1&x_3&&&r_3\\ &&1&&& \\
&&&1&y_1&y_2\\ &&&&1&y_2\\ &&&&&1 \end{pmatrix} \begin{pmatrix}
1&&&&&\\ &1&&&&\\ &&1&&&\\ &z_1&z_2&1&& \\ &&z_3&&1&\\&&&&&1
\end{pmatrix} \right )\psi(x_1+x_3+y_1+y_3)d(...)
\end{equation}
is zero for all choice of data. Expand \eqref{speh3} along the
unipotent group consisting of all matrices
$I_6+r_4e_{1,4}+r_5e_{2,5}+r_6e_{3,6}$ where $r_i\in F\backslash
{\bf A}$. Conjugating by a suitable discrete matrix and collapsing
summations with integrations, \eqref{speh3} is equal to
\begin{equation}\label{speh4}
\int E_\tau\left (
\begin{pmatrix} 1&x_1&x_2&r_4&r_1&r_2\\ &1&x_3&&r_5&r_3\\ &&1&&&r_6 \\
&&&1&y_1&y_2\\ &&&&1&y_2\\ &&&&&1 \end{pmatrix} \begin{pmatrix}
1&&&&&\\ &1&&&&\\ &&1&&&\\ &z_1&z_2&1&& \\ &&z_3&&1&\\&&&&&1
\end{pmatrix} \right )\psi(x_1+x_3+y_1+y_3)d(...)
\end{equation}
where the variables $z_i$ are integrated over ${\bf A}$ and all
the rest over $F\backslash {\bf A}$. Thus we conclude that
integral \eqref{speh4} is zero for all choice of data. As in
\cite{Ga-S} section 7, it follows that we may ignore the
integration over the $z_i$ variables. That is, the inner
integration is zero for all choice of data. Expand the inner
integration in \eqref{speh4} along $I_6+r_7e_{2,4}+r_8e_{3,5}$
where $r_i\in F\backslash {\bf A}$. We claim that except the
constant term, all other terms contribute zero to the expansion.
Indeed, this follows from the fact that each other term produces a
Fourier coefficient for $E_\tau$ which corresponds to a unipotent
orbit which is greater than or not related to $(3^2)$. Next,
expanding along $I+r_9e_{3,4}$ with $r_9\in F\backslash {\bf A}$,
all terms, except the constant term, contribute zero. Thus we may
deduce that the integral
\begin{equation}\label{speh5}
\int\limits_{(F\backslash {\bf A})^{12}}E_\tau\left (
\begin{pmatrix} 1&x_1&x_2&r_4&r_1&r_2\\ &1&x_3&r_7&r_5&r_3\\ &&1&r_9&r_8&r_6 \\
&&&1&y_1&y_2\\ &&&&1&y_2\\ &&&&&1
\end{pmatrix}\right )\psi(x_1+x_3+y_1+y_3)d(...)
\end{equation}
is zero for all choice of data. However, from the definition of
the Eisenstein series $E_\tau(g,s)$, see beginning of section 5,
it follows that this last integral cannot be zero for all choice
of data. Thus we obtained a contradiction.    \hfill
$\blacksquare$

\end{document}